\ifdefined\WettingReadingVersion
\documentclass[hidelinks,onefignum,onetabnum,final]{siamart251216}
\else
\documentclass[hidelinks,onefignum,onetabnum,final]{siamart251216}
\fi
\usepackage{amssymb,mathrsfs,mathtools}
\usepackage{booktabs}
\usepackage{amsmath,amssymb}
\usepackage{graphicx}
\usepackage{bm,stmaryrd}
\usepackage{enumitem}
\usepackage{placeins}

\newsiamremark{remark}{Remark}
\newsiamremark{hypothesis}{Hypothesis}
\crefname{hypothesis}{Hypothesis}{Hypotheses}
\newsiamremark{assumption}{Assumption}
\crefname{assumption}{Assumption}{Assumptions}
\newsiamthm{claim}{Claim}
\newsiamremark{fact}{Fact}
\crefname{fact}{Fact}{Facts}

\usepackage{xcolor}

\newcommand{\bn}{\mathbf{n}}

\newcommand{\cE}{\mathcal{E}}

\newcommand{\cV}{\mathcal{V}}

\newcommand{\eps}{\varepsilon}

\renewcommand{\d}{\mathrm{d}}
\newcommand{\dx}{\,\d x}

\newcommand{\dS}{\,\mathrm{d} S}

\newcommand{\ptial}[1]{ \partial_{#1} }
\newcommand{\pt}{\ptial{t}}

\newsiamthm{scheme}{Scheme}

\headers{Nonlocal Allen--Cahn Wetting}{Y. Liu, N. Lu, Z. Qiao, and Y. Zhang}

\title{A Nonlocal Allen--Cahn Model for Volume-Constrained Wetting and Its Multi-Structure-Preserving Discretization\thanks{
	Funding: Y. Liu is partially supported by the NSFC No.~12601426, the NSF of Jiangsu Province of China under Grant No.~BK20240572, and the China Postdoctoral Science Foundation under Grant No.~2025M773078. N. Lu is partially supported by the NSFC  No. 12601780. Z. Qiao's work is partially supported by NSFC/RGC Joint Research Scheme N\_PolyU5145/24 and the Hong Kong Polytechnic University RISE project 4-CDDQ. Y. Zhang is partially supported by the NSFC No.~12301554.
}
}

\author{
Yadong Liu\thanks{School of Mathematical Sciences, and Ministry of Education Key Laboratory of NSLSCS, and Key Laboratory of Jiangsu Provincial Universities of FDMTA, Nanjing Normal University, Nanjing 210023, P. R. China
	(\email{ydliu@njnu.edu.cn})}
\and Nan Lu\thanks{Department of Mathematics, Shanghai Normal University, Shanghai 200234, P. R. China
	(\email{lunan@shnu.edu.cn})}
\and Zhonghua Qiao\thanks{Department of Applied Mathematics and Research Institute for Smart Energy, the Hong Kong Polytechnic University, Hung Hom, Kowloon, Hong Kong (\email{zhonghua.qiao@polyu.edu.hk})}
\and Yuze Zhang\thanks{School of Mathematical Sciences, and Ministry of Education Key Laboratory of NSLSCS, and Key Laboratory of Jiangsu Provincial Universities of FDMTA, Nanjing Normal University, Nanjing 210023, P. R. China
	(\email{05439@njnu.edu.cn}, corresponding author)}
}

\date{\today}

\begin{document}

\maketitle

\begin{abstract}
At finite interface thickness, small droplets can shrink or disappear during wetting relaxation even when the phase-field mass is conserved. We develop a nonlocal Allen--Cahn model that couples a nonlinear diffuse-volume constraint to dynamic wall relaxation through Onsager's principle. The constraint gradient vanishes at the pure-phase values and does zero work in the energy balance, preserving the original bulk and wall dissipation structure. We construct a Compatible Variational Chain-Rule (CVC) scheme using the same discrete derivative of the volume density in the correction and in the volume increment. The scheme exactly conserves the nonlinear diffuse volume and reproduces the bulk--wall dissipation identity for the original discrete energy. Under a coercivity condition, the energy also yields a mesh-dependent maximum-norm bound. Numerical tests show second-order convergence in space and time. In flat-wall tests, the proposed model improves small-droplet retention and geometric-volume preservation and reproduces equilibrium cap shapes and contact angles. Additional tests illustrate wetting on chemically heterogeneous and rough substrates.
\end{abstract}



\section{Introduction}\label{sec:introduction}

For a nonvolatile, incompressible droplet relaxing on a solid wall, the shape and contact-line positions change while the volume should remain fixed \cite{LuXu2021}. A mass constraint on the phase variable does not necessarily enforce this requirement at finite interface thickness: droplet shrinkage can be compensated by small shifts of the phase variable over the surrounding region. Yue, Zhou, and Feng analyzed this mechanism for Cahn--Hilliard dynamics \cite{YueZhouFeng2007}. For the nonlocal Allen--Cahn equation, dissolution of small features
into the surrounding phase has been linked to the spatially uniform mass correction~\cite{KimLeeChoi2014}. In our flat-wall tests, a small droplet disappears during wetting relaxation while the global phase-field mass is conserved to solver accuracy. We seek a correction that limits such loss of geometric volume and retains the original bulk--wall energy-dissipation balance.

Volume errors affect comparisons of droplet shapes and contact angles. Surface texture produces distinct wetting states \cite{LafumaQuere2003}; control of wettability enables self-cleaning \cite{Blossey2003}, and wettability gradients drive droplet motion \cite{ChaudhuryWhitesides1992}. Their characterization requires distinguishing equilibrium, advancing, and receding contact angles \cite{HuhtamakiEtAl2018}. Variational contact-line models and analyses of hysteresis describe the wall interactions involved \cite{QianWangSheng2006,XuZhaoWang2019}.
The classical mass-conserving Allen--Cahn model preserves the linear integral of the phase field, commonly referred to as ``mass'', and dissipates the bulk energy \cite{RS1992}. Its spatially uniform correction, however, can shift the bulk values and thereby compensate for changes in droplet volume. Brassel and Bretin \cite{BresselBretin2011} introduced the weight $\sqrt{2W(u)}$, with $W$ the bulk potential, to suppress forcing of the pure phases and improve volume preservation. Since this weight vanishes at $u=\pm1$, the correction acts primarily within the transition layer. Their model nevertheless keeps the linear mass constraint, whose $L^2$ gradient is constant. The weighted correction therefore does not match the gradient of the conserved functional, so the standard dissipation identity for the original energy does not generally follow \cite{AA2014}.

	To combine this interface weight with the original energy-dissipation structure, we choose a conserved functional whose variational derivative coincides with $\sqrt{2W(u)}$. Following Golovaty \cite{Golovaty1997}, we impose the nonlinear constraint
	\begin{align}\label{eq:diffuse-volume}
		H(s)=\int_0^s\sqrt{2W(r)}\,\mathrm{d}r,
		\qquad
		\mathcal V(u)=\int_\Omega H(u)\,\mathrm{d}x,
	\end{align}
	where $W$ is a nonnegative, balanced double-well potential satisfying $W(\pm1)=0$ and $W>0$ on $(-1,1)$. The correction $\lambda(t)H'(u)=\lambda(t)\sqrt{2W(u)}$ now acts along the gradient of $\mathcal V$. Note that by the conservation of $\mathcal V$, it follows $\int_\Omega H'(u)\partial_tu\,\mathrm{d}x=0$. Thus, essentially the new correction makes \emph{zero contribution} to the energy balance. This choice retains the vanishing of the correction in the pure phases while ensuring a constrained-gradient-flow structure for the original energy \cite{KL2025,Poiatti2026,Takasao2023}. The precise evolution and multiplier are given in \eqref{eq:constrained-allen-cahn} and \eqref{eq:lambda}, and the resulting energy-dissipation identity is stated in \eqref{eq:energy-dissipation}.
		The corresponding normalized diffuse-volume is
	\begin{align}\label{eq:intro-normalized-volume}
		V_H(u)
		=\int_\Omega\frac{H(u)-H(-1)}{H(1)-H(-1)}\,\mathrm{d}x
		=\frac{\mathcal V(u)-H(-1)|\Omega|}{H(1)-H(-1)}.
	\end{align}
	On the fixed domain $\Omega$, its conservation is equivalent to that of $\mathcal V$, and $V_H(2\mathbf1_D-1)=|D|$. Numerical comparisons in Figure~\ref{fig:reference-model-comparison} demonstrate improved droplet retention, although exact conservation concerns $V_H$, rather than the geometric volume enclosed by a level set at finite $\varepsilon$. Note that the choice $H'=\sqrt{2W}$ also makes this density compatible with the prescribed interfacial energy. Further geometric and energetic interpretations of the nonlinear volume constraint are given in Appendix~\ref{app:diffuse-volume}.

We extend this formulation to wetting by incorporating wall energy and finite-rate boundary relaxation through Onsager's principle \cite{QianWangSheng2006,LuXu2021,ChenLiuXu2025}. For the quartic potential used here, $g(x,s)=-\cos\theta(x)H(s)$ on $[-1,1]$, so the volume constraint and wall energy use the same nonlinear function. The resulting model couples contact-line relaxation to conservation of $\mathcal V$ and dissipation of the original bulk--wall energy.

Numerical schemes for phase-field models with moving contact lines have been developed using stabilized and convex-splitting approaches, mixed variational formulations, and treatments of dynamic boundary relaxation \cite{ShenYang2015,LiQiaoTang2016,GomezHughes2011,ShenYangYu2015,YangYu2018}. Auxiliary-variable methods, including IEQ and SAV, provide further approaches to energy-stable time discretization \cite{ShenXuYang2019,AkrivisLiLi2019}. Methods designed for global constraints and for dissipation of the original energy address additional aspects of conservation and energy stability \cite{ChengShen2020,ChengLiuShen2020}.
Existing numerical methods for phase-field models with moving contact lines have addressed energy stability, global constraints, and the treatment of dynamic boundary conditions through stabilized, convex-splitting, auxiliary-variable, and multiplier approaches \cite{ShenYang2015,LiQiaoTang2016,GomezHughes2011,ShenYangYu2015,YangYu2018,ShenXuYang2019,ChengShen2020,ChengLiuShen2020}. Discrete variational derivatives and discrete-gradient methods provide another approach by enforcing exact discrete chain rules \cite{Furihata2001,CelledoniEtAl2012}. Related fully discrete methods have also combined conservation and energy dissipation for other classes of evolution equations \cite{HuVanFleetWan2025}. These developments motivate a more specific question for the present dynamic-boundary model: can the nonlinear diffuse-volume and the original bulk--wall energy-dissipation rate be preserved simultaneously at the fully discrete level? This is more restrictive than imposing a linear mass constraint or establishing energy decay alone. The volume constraint is defined by a nonlinear functional of the phase field, while the energy balance couples the bulk evolution to the dynamic boundary relaxation. Consequently, the constraint work and the boundary contribution must be handled consistently within the same discretization.

To this end, we develop a Compatible Variational Chain-Rule (CVC) scheme by combining summation-by-parts spatial operators with discrete derivatives satisfying exact chain rules for the bulk potential, the nonlinear diffuse-volume density, and the wall energy. The same discrete constraint that ensures volume conservation cancels the multiplier contribution to the energy balance, while the compatible boundary discretization retains the wall dissipation. The main contributions are as follows.

\begin{enumerate}
	
	\item \textbf{Compatible wetting modeling.}
	An interface-weighted nonlinear volume constraint is coupled
	to dynamic wall relaxation through Onsager's principle,
	reducing sensitivity to bulk shifts while retaining the
	original energy-dissipation law.
	
	\item \textbf{Simultaneous structure preservation.}
	Summation-by-parts operators and exact discrete chain rules
	for the bulk potential, volume density, and wall energy yield
	a scheme that conserves $\mathcal V_h$ and satisfies the full
	discrete dissipation identity for the original energy, including
	bulk and boundary terms. We establish formal second-order
	temporal consistency on smooth branches and a mesh-dependent
	maximum-norm estimate \cite{DuJuLiQiao2019,DuJuLiQiao2021} under a stated coercivity condition.
	
	\item \textbf{Verification and wetting accuracy.}
	Experiments verify second-order convergence and the structural
	identities, demonstrate improved geometric-volume preservation
	and small-droplet retention over the mass-conserving
	model, and assess equilibrium shapes and contact angles.
	Further tests cover pinning, chemically directed migration
	and coalescence, and rough-surface wetting.

\end{enumerate}

Section~\ref{sec:model} presents the model and continuous properties; Section~\ref{sec:discretization} develops the discretization and estimates. Numerical experiments and conclusions follow in Sections~\ref{sec:numerical} and~\ref{sec:conclusion}, respectively.

\section{Nonlocal Allen--Cahn wetting model and continuous structure}\label{sec:model}

\subsection{Rayleighian and constitutive choices}

Let $\Omega$ be the fluid domain and let $\Gamma\subseteq\partial\Omega$ be the solid wall on which wetting relaxation occurs. The remaining boundaries are periodic or carry homogeneous Neumann conditions, so that they contribute no boundary energy or dissipation. We use the double-well potential and wall energy density
\begin{align*}
	W(s)=\frac14(s^2-1)^2,\qquad
	g(x,s)=-\frac{\sqrt2}{6}\cos\theta(x)\,(3s-s^3).
\end{align*}
Here $\theta(x)$ is the prescribed, time-independent Young angle on the wall. We write $g(u)$ and $g'(u)$ for $g(x,u(x))$ and $\partial_s g(x,u(x))$, respectively. The functions $H$ and $\cV$ are defined in \eqref{eq:diffuse-volume}, with $H'=\sqrt{2W}$.

At a fixed sufficiently regular state $u$, let $v$ denote a trial rate, with its trace on $\Gamma$ describing the wall relaxation, and let $\lambda\in\mathbb R$ be a spatially constant Lagrange multiplier. We define the constrained \emph{Rayleighian} by
\begin{align}
	\mathscr R_\eps(u;v,\lambda)
	&\coloneqq
	\mathcal{D}_\eps(v)
	+
	D\cE_\eps(u)[v]
	-
	\eps\lambda D\cV(u)[v]
	\nonumber\\
	&=
	\frac{\eps}{2}\int_\Omega |v|^2\dx
	+\frac12\int_\Gamma |v|^2\dS
	+
	\int_\Omega
	\left(
	\eps\nabla u\cdot\nabla v
	+\frac{1}{\eps}W'(u)v
	\right)\dx
	\nonumber\\
	&\quad
	+
	\int_\Gamma g'(u)v\dS
	-
	\eps\lambda\int_\Omega \sqrt{2W(u)}v\dx.
	\label{eq:rayleighian}
\end{align}
Here the energy and dissipation potential are respectively
\begin{align}
	\cE_\eps(u)
	&\coloneqq
	\int_\Omega
	\left(\frac{\eps}{2}|\nabla u|^2+\frac1\eps W(u)\right)\dx
	+\int_\Gamma g(u)\dS,
	\label{eq:energy}\\
	\mathcal D_\eps(v)
	&\coloneqq
	\frac{\eps}{2}\int_\Omega |v|^2\dx
	+\frac12\int_\Gamma |v|^2\dS,
	\label{eq:dissipation-potential}
\end{align}
with directional derivatives of the energy and volume satisfying
\begin{align*}
	D\cE_\eps(u)[v]
	&=\int_\Omega
	\left(\eps\nabla u\cdot\nabla v+\frac1\eps W'(u)v\right)\dx
	+\int_\Gamma g'(u)v\dS,\\
	D\cV(u)[v]
	&=\int_\Omega H'(u)v\dx
	=\int_\Omega\sqrt{2W(u)}v\dx.
\end{align*}
Onsager's principle determines the physical rate by stationarity of \eqref{eq:rayleighian} with respect to $v$ and $\lambda$. Equivalently, it minimizes $\mathcal D_\eps(v)+D\cE_\eps(u)[v]$ over rates satisfying $D\cV(u)[v]=0$, with $\lambda$ enforcing this constraint. We now derive the evolution equations, volume conservation, and energy dissipation from these stationarity conditions.

\subsection{Evolution equations and continuous laws from stationarity}

Let $\xi\in H^1(\Omega)$ be an arbitrary variation of the rate $v$, with the same periodicity where applicable. Stationarity with respect to $v$ gives
\begin{align*}
	0
	&=\left.\frac{\mathrm d}{\mathrm ds}
	\mathscr R_\eps(u;v+s\xi,\lambda)\right|_{s=0}\\
	&=\int_\Omega
	\left(
	\eps v\xi+\eps\nabla u\cdot\nabla\xi
	+\frac1\eps W'(u)\xi
	-\eps\lambda\sqrt{2W(u)}\xi
	\right)\dx
	+\int_\Gamma\bigl(v+g'(u)\bigr)\xi\dS.
\end{align*}
Integrating the gradient term by parts yields
\begin{align*}
	0
	&=\eps\int_\Omega
	\left(
	v-\Delta u+\frac1{\eps^2}W'(u)
	-\lambda\sqrt{2W(u)}
	\right)\xi\dx
	+\int_\Gamma
	\left(v+\eps\partial_{\bn}u+g'(u)\right)\xi\dS.
\end{align*}
Choosing $\xi$ with compact support in $\Omega$, and allowing arbitrary traces on $\Gamma$, gives
\begin{align*}
	v&=\Delta u-\frac1{\eps^2}W'(u)
	+\lambda\sqrt{2W(u)}
	&&\text{in }\Omega,\\
	v&=-\bigl(\eps\partial_{\bn}u+g'(u)\bigr)
	&&\text{on }\Gamma.
\end{align*}
The dynamic boundary condition follows from the boundary dissipation in \eqref{eq:rayleighian}; omitting that quadratic boundary term would give the static natural condition $\eps\partial_{\bn}u+g'(u)=0$. Taking $v=\pt u_\eps$ and $\lambda=\lambda_\eps(t)$, we obtain
\begin{align}
	\left\{
	\begin{aligned}
		\pt u_\eps
		&=\Delta u_\eps-\frac1{\eps^2}W'(u_\eps)
		+\lambda_\eps\sqrt{2W(u_\eps)}
		&&\text{in }\Omega\times(0,T),\\
		\pt u_\eps
		&=-\bigl(\eps\partial_{\bn}u_\eps+g'(u_\eps)\bigr)
		&&\text{on }\Gamma\times(0,T).
	\end{aligned}
	\right.
	\label{eq:constrained-allen-cahn}
\end{align}

Next, stationarity with respect to $\lambda$ gives
\begin{align*}
	0=\partial_\lambda\mathscr R_\eps(u;v,\lambda)
	=-\eps D\cV(u)[v]
	=-\eps\int_\Omega\sqrt{2W(u)}v\dx.
\end{align*}
Along the evolution, this becomes
\begin{align}
	\int_\Omega\sqrt{2W(u_\eps)}\pt u_\eps\dx=0.
	\label{eq:lambda-constraint}
\end{align}
Since $H'=\sqrt{2W}$, the chain rule yields
\begin{align*}
	\frac{\mathrm d}{\mathrm dt}\cV(u_\eps(t))
	=\int_\Omega H'(u_\eps)\pt u_\eps\dx=0,
\end{align*}
and hence
\begin{align}
	\cV(u_\eps(t))=\cV(u_\eps(0)).
	\label{eq:volume-conservation}
\end{align}
After a fixed shift and positive scaling, $\cV$ approximates the liquid volume as explained in \eqref{eq:intro-normalized-volume}. At finite interface thickness, \eqref{eq:volume-conservation} does not imply exact preservation of the volume enclosed by the zero level set.

The multiplier is determined by the equation. Namely, by substituting the bulk equation into \eqref{eq:lambda-constraint}, we have
\begin{align*}
	0
	=\int_\Omega
	\left(\Delta u_\eps-\frac1{\eps^2}W'(u_\eps)\right)
	\sqrt{2W(u_\eps)}\dx
	+\lambda_\eps\int_\Omega 2W(u_\eps)\dx,
\end{align*}
which then implies
\begin{align}
	\lambda_\eps(t)
	=-\frac{
		\int_\Omega
		\left(\Delta u_\eps-\frac1{\eps^2}W'(u_\eps)\right)
		\sqrt{2W(u_\eps)}\dx
	}{
		\int_\Omega 2W(u_\eps)\dx
	},
	\label{eq:lambda}
\end{align}
provided that $\int_\Omega 2W(u_\eps)\dx>0$. 

Finally, choosing the rate itself as the variation, $\xi=v$, in the stationarity condition, one obtains
\begin{align*}
	0=2\mathcal D_\eps(v)+D\cE_\eps(u)[v]
	-\eps\lambda D\cV(u)[v].
\end{align*}
The multiplier contribution vanishes by stationarity with respect to $\lambda$. Setting $v=\pt u_\eps$ and applying the chain rule therefore yields
\begin{align}
	\frac{\mathrm d}{\mathrm dt}\cE_\eps(u_\eps)
	=-2\mathcal D_\eps(\pt u_\eps)=-\eps\|\pt u_\eps\|_{L^2(\Omega)}^2
	-\|\pt u_\eps\|_{L^2(\Gamma)}^2.\label{eq:energy-dissipation}
\end{align}
Thus both diffuse-volume conservation and dissipation of the original bulk--wall energy follow directly from the same Rayleighian, and \eqref{eq:constrained-allen-cahn}--\eqref{eq:lambda} is the corresponding constrained gradient flow. The derivation applies to sufficiently regular evolutions for which the trace identities, chain rules, and integrations by parts are valid.

\subsection{Conservation, dissipation, and the physical interval}

The variational derivation is formal unless the state and its trace have the regularity needed for the integrations by parts. The following result states the continuous properties for classical solutions; it is not an existence theorem for the nonlocal problem.

\begin{theorem}[Continuous structure]\label{thm:continuous-structure}
	Let $u$ be a classical solution of \eqref{eq:constrained-allen-cahn}--\eqref{eq:lambda} on $[0,T]$, with the boundary conditions specified above. Assume that the integrations by parts below are valid, $\int_\Omega 2W(u(t))\,\mathrm dx>0$ for every $t\in[0,T]$, and $\lambda_\eps\in L^\infty(0,T)$. The prescribed wall energy is independent of time. Then
	\begin{equation}
		\cV(u(t))=\cV(u(0)),\qquad
		\cE_\eps(u(t))+
		\int_0^t\left(\eps\|u_t\|_{L^2(\Omega)}^2+
		\|u_t\|_{L^2(\Gamma)}^2\right)\,\mathrm ds
		=\cE_\eps(u(0)).
		\label{eq:continuous-integrated-laws}
	\end{equation}
	If $-1\le u(0)\le1$ in $\overline\Omega$, then $-1\le u(t)\le1$ in $\overline\Omega$ for $0\le t\le T$.
\end{theorem}

\begin{proof}
	Set $q=H'(u)$ and $A=\Delta u-\eps^{-2}W'(u)$. The multiplier satisfies $\lambda_\eps=-(A,q)_{L^2(\Omega)}/\|q\|_{L^2(\Omega)}^2$; hence
	\[
	\frac{\mathrm d}{\mathrm dt}\cV(u)
	=(q,u_t)_{L^2(\Omega)}=(q,A+\lambda_\eps q)_{L^2(\Omega)}=0.
	\]
	Differentiating the total energy and integrating by parts gives
	\begin{align*}
		\frac{\mathrm d}{\mathrm dt}\cE_\eps(u)
		&=\int_\Omega(-\eps\Delta u+\eps^{-1}W'(u))u_t\,\mathrm dx
		+\int_\Gamma(\eps\partial_{\bn}u+g'(u))u_t\,\mathrm dS\\
		&=-\eps\|u_t\|_{L^2(\Omega)}^2-\|u_t\|_{L^2(\Gamma)}^2
		+\eps\lambda_\eps(q,u_t)_{L^2(\Omega)}.
	\end{align*}
	The last term vanishes by the constraint. Integration in time proves \eqref{eq:continuous-integrated-laws}.
	
	For the upper bound, define $w=(u-1)_+$ and $f(t,s)=-\eps^{-2}W'(s)+\lambda_\eps(t)H'(s)$. The quartic potential and the chosen wall energy satisfy $f(t,1)=g'(1)=0$. Since $u$ is bounded on the classical existence interval and $\lambda_\eps$ is bounded, $f(t,\cdot)$ and $g'$ are uniformly Lipschitz on the range needed in this argument. Testing the bulk equation by $\eps w$ and the boundary equation by its trace $w$ cancels the normal-flux terms:
	\begin{align*}
		\frac12\frac{\mathrm d}{\mathrm dt}
		\Big(\eps\|w\|_{L^2(\Omega)}^2&+\|w\|_{L^2(\Gamma)}^2\Big)
		=-\eps\|\nabla w\|_{L^2(\Omega)}^2
		+\eps\int_\Omega(f(t,u)-f(t,1))w\,\mathrm dx\\
		&\quad-\int_\Gamma(g'(u)-g'(1))w\,\mathrm dS\leq C\left(\eps\|w\|_{L^2(\Omega)}^2+
		\|w\|_{L^2(\Gamma)}^2\right).
	\end{align*}
	The same estimate holds pointwise in the wall coordinate when $g$ depends on position through a bounded prescribed $\cos\theta$. Since $w(0)=0$, Gronwall's inequality yields $w=0$. Testing instead with the negative of $z=(-1-u)_+$, and using $f(t,-1)=g'(-1)=0$, gives the lower bound in the same way.
\end{proof}


\section{Compatible discretization and discrete structure}\label{sec:discretization}
Motivated by the constrained variational structure of the proposed model, we construct a second-order fully discrete Compatible Variational Chain-Rule (CVC) scheme. The spatial discretization preserves the bulk--boundary variational coupling through a discrete summation-by-parts identity. For the temporal discretization, we employ discrete derivatives satisfying exact chain-rule identities that are compatible with the variational relations of the semi-discrete system. Consequently, the fully discrete scheme exactly conserves the discrete diffuse volume and satisfies an exact dissipation law for the discrete counterpart of the original energy.
\subsection{Spatial operators and summation by parts}

Let $\Omega=[a,b]\times[c,d]$ be a rectangular domain, where the dynamic boundary condition is imposed on the lower boundary $\Gamma=\{(x,c):a\le x\le b\}$. Periodic boundary conditions are imposed in the $x$-direction, while the homogeneous Neumann boundary condition is imposed on the upper boundary.

For two positive integers $N_x$ and $N_y$, the domain $\Omega$ is uniformly partitioned with mesh sizes $\Delta x=(b-a)/N_x$, $\Delta y=(d-c)/N_y$. The spatial grid points are denoted by $\Omega_h = \left\{ (x_j,y_k)\,\middle|\, x_j=a+j\Delta x,\; y_k=c+k\Delta y,\; 0\le j\le N_x-1,\; 0\le k\le N_y \right\}$. For a given grid function $f_{j,k}$, we introduce the following difference operators:
\begin{align*}
	&\delta_x^+ f_{j,k}
	=
	\frac{f_{j+1,k}-f_{j,k}}{\Delta x},
	\qquad
	&&\delta_y^+ f_{j,k}
	=
	\frac{f_{j,k+1}-f_{j,k}}{\Delta y},
	\\
	&\delta_x^{\langle2\rangle} f_{j,k}
	=
	\frac{f_{j+1,k}-2f_{j,k}+f_{j-1,k}}{(\Delta x)^2},
	\qquad
	&&\delta_y^{\langle2\rangle} f_{j,k}
	=
	\frac{f_{j,k+1}-2f_{j,k}+f_{j,k-1}}{(\Delta y)^2},
	\\
	&\delta_y^{\langle1\rangle} f_{j,k}
	=
	\frac{f_{j,k+1}-f_{j,k-1}}{2\Delta y}.
\end{align*}
The periodicity in the $x$-direction and the homogeneous Neumann condition on the upper boundary are imposed through
\begin{equation*}
	f_{-1,k}=f_{N_x-1,k},
	\qquad
	f_{N_x,k}=f_{0,k},
	\qquad
	f_{j,N_y+1}=f_{j,N_y-1},
\end{equation*}
while a ghost value $f_{j,-1}$ is introduced to accommodate the centered difference approximation on the lower boundary. The discrete Laplace operator and the outward normal derivative on $\Gamma$ are then defined by
\begin{equation*}
	\Delta_h f_{j,k}
	=
	\delta_x^{\langle2\rangle}f_{j,k}
	+
	\delta_y^{\langle2\rangle}f_{j,k},\qquad
	(\partial_{\bn,h}f)_{j,0}
	=
	-\delta_y^{\langle1\rangle}f_{j,0}.
\end{equation*}
The corresponding discrete inner products and norms in $\Omega$ and on $\Gamma$, together with the discrete gradient inner product, are defined as
\begin{align*}
	&\langle f,g\rangle_{\Omega,h}
	=
	\Delta x\Delta y
	\sum_{j=0}^{N_x-1}
	\sum_{k=0}^{N_y}
	\omega_k f_{j,k}g_{j,k},
	\qquad
	\|f\|_{\Omega,h}
	=
	\sqrt{\langle f,f\rangle_{\Omega,h}},
	\\
	&\langle f,g\rangle_{\Gamma,h}
	=
	\Delta x
	\sum_{j=0}^{N_x-1}
	f_{j,0}g_{j,0},
	\qquad
	\|f\|_{\Gamma,h}
	=
	\sqrt{\langle f,f\rangle_{\Gamma,h}},
	\\
	&\langle \nabla_h f,\nabla_h g\rangle_h
	=
	\Delta x\Delta y
	\sum_{j=0}^{N_x-1}
	\sum_{k=0}^{N_y}
	\omega_k
	(\delta_x^+f_{j,k})(\delta_x^+g_{j,k})
	\\
	&\hspace{2em}
	+
	\Delta x\Delta y
	\sum_{j=0}^{N_x-1}
	\sum_{k=0}^{N_y-1}
	(\delta_y^+f_{j,k})(\delta_y^+g_{j,k}),
\end{align*}
where the weights $\omega_k$ are defined by $\omega_k=
\begin{cases}
	1/2, & k=0,N_y,\\
	1, & 1\le k\le N_y-1,
\end{cases}$
corresponding to the trapezoidal rule in the $y$-direction. For grid functions satisfying the above boundary conventions, the following summation-by-parts identity holds:
\begin{equation}
	\langle \Delta_h f,g\rangle_{\Omega,h}
	=
	-\langle \nabla_h f,\nabla_h g\rangle_h
	+
	\langle \partial_{\bn,h}f,g\rangle_{\Gamma,h}.
	\label{eq:discrete-sbp}
\end{equation}
Indeed, periodic summation in the $x$-direction gives
\[
\Delta x\sum_{j=0}^{N_x-1}
(\delta_x^{\langle2\rangle}f)_{j,k}g_{j,k}
=
-\Delta x\sum_{j=0}^{N_x-1}
(\delta_x^+f)_{j,k}(\delta_x^+g)_{j,k}.
\]
In the $y$-direction, summation by parts yields
\begin{align*}
	&\Delta y\sum_{k=0}^{N_y}\omega_k
	(\delta_y^{\langle2\rangle}f)_{j,k}g_{j,k}
	\\
	={}&
	-\Delta y\sum_{k=0}^{N_y-1}
	(\delta_y^+f)_{j,k}(\delta_y^+g)_{j,k}
	+\frac{f_{j,N_y+1}-f_{j,N_y-1}}{2\Delta y}g_{j,N_y}
	-\frac{f_{j,1}-f_{j,-1}}{2\Delta y}g_{j,0}.
\end{align*}
The upper boundary term vanishes by the homogeneous Neumann condition, while the lower boundary term equals $(\partial_{\bn,h}f)_{j,0}g_{j,0}$. Multiplying by $\Delta x$ and summing over $j$ gives \eqref{eq:discrete-sbp}. 
The identity also holds after eliminating the lower
ghost values through the dynamic boundary condition. For sufficiently smooth solutions, the ghost-point
construction yields second-order-consistent boundary
approximations while preserving the discrete
summation-by-parts identity \eqref{eq:discrete-sbp}.

For $0\le j<N_x$ and $0\le k\le N_y$, the spatially semi-discrete bulk equation and dynamic boundary condition are given by
\begin{equation}\label{eq:semi-discrete}
	\left\{
	\begin{aligned}
		&\frac{\mathrm{d}U_{j,k}}{\mathrm{d}t}
		=
		\Delta_h U_{j,k}
		-\frac{1}{\eps^2}W'(U_{j,k})
		+\lambda_{\eps,h}\sqrt{2W(U_{j,k})},
		\\
		&\frac{\mathrm{d}U_{j,0}}{\mathrm{d}t}
		=
		-\left(
		\eps(\partial_{\bn,h}U)_{j,0}
		+g'(U_{j,0})
		\right).
	\end{aligned}
	\right.
\end{equation}
where the discrete Lagrange multiplier is determined by
\begin{equation}
	\lambda_{\eps,h}
	=
	\frac{
		\left\langle
		\frac{1}{\eps^2}W'(U)-\Delta_h U,
		\sqrt{2W(U)}
		\right\rangle_{\Omega,h}
	}{
		\left\|\sqrt{2W(U)}\right\|_{\Omega,h}^2
	}.
	\label{eq:semi-discrete-lambda}
\end{equation}
The lower ghost values are determined by combining the bulk equation at $k=0$ with the dynamic boundary condition. The elimination of these values in the fully discrete scheme is described in Remark~\ref{rem:ghost-elimination}. The semi-discrete system \eqref{eq:semi-discrete} satisfies the following discrete volume conservation and energy dissipation laws.
\begin{theorem}
	\label{thm:semi-discrete-structure}
	The semi-discrete system \eqref{eq:semi-discrete} satisfies
	\begin{equation*}
		\frac{\mathrm{d}}{\mathrm{d}t}\mathcal{V}_h(U)=0,
		\qquad
		\frac{\mathrm{d}}{\mathrm{d}t}\mathcal{E}_{\eps,h}(U)
		=
		-\eps\left\|\frac{\mathrm{d}U}{\mathrm{d}t}\right\|_{\Omega,h}^2
		-
		\left\|\frac{\mathrm{d}U}{\mathrm{d}t}\right\|_{\Gamma,h}^2
		\le 0,
	\end{equation*}
	where $\mathcal{V}_h$ and $\mathcal{E}_{\eps,h}$ are the discrete counterparts of the diffuse-volume $\mathcal{V}$ and the original energy $\mathcal{E}_\eps$, respectively, defined by
	\begin{align}
		&\mathcal{V}_h(U)
		=
		\langle H(U),1\rangle_{\Omega,h},
		\notag
		\\
		&\mathcal{E}_{\eps,h}(U)
		=
		\frac{\eps}{2}\|\nabla_h U\|_h^2
		+
		\frac{1}{\eps}\langle W(U),1\rangle_{\Omega,h}
		+
		\langle g(U),1\rangle_{\Gamma,h}.
		\label{eq:discrete-energy}
	\end{align}
\end{theorem}
\begin{proof}
	The proof follows the same argument as that of Theorem~\ref{thm:continuous-structure}, replacing the integration-by-parts formula by the discrete summation-by-parts identity \eqref{eq:discrete-sbp}.
\end{proof}
\subsection{Compatible variational chain rules and time stepping}
	Discrete variational derivatives have been used to construct structure-preserving finite-difference schemes \cite{Furihata2001}, and related chain-rule constructions arise in the average-vector-field approach to dissipative PDEs \cite{CelledoniEtAl2012}. The present construction is tailored to the constrained bulk--boundary variational structure of the proposed model. Let $t_n=n\tau$, $n=0,\ldots,N$, with $\tau=T/N$, and define
	\begin{equation*}
		U^{n+\frac12}
		=\frac{U^{n+1}+U^n}{2},
		\qquad
		\delta_t U^{n+\frac12}
		=\frac{U^{n+1}-U^n}{\tau}.
	\end{equation*}
	Along a trajectory of the semi-discrete system, the
discrete volume and energy satisfy the following
chain-rule identities:
	\begin{align*}
		&\frac{\mathrm d}{\mathrm dt}\mathcal V_h(U)
		=
		\left\langle
		H'(U),\frac{\mathrm dU}{\mathrm dt}
		\right\rangle_{\Omega,h},
		\\
		&\frac{\mathrm d}{\mathrm dt}\mathcal E_{\eps,h}(U)
		=
		\eps\left\langle
		\nabla_h U,\nabla_h\frac{\mathrm dU}{\mathrm dt}
		\right\rangle_h+
		\frac{1}{\eps}
		\left\langle
		W'(U),\frac{\mathrm dU}{\mathrm dt}
		\right\rangle_{\Omega,h}
		+
		\left\langle
		g'(U),\frac{\mathrm dU}{\mathrm dt}
		\right\rangle_{\Gamma,h}.
	\end{align*}
	To extend these variational identities to the fully
discrete level, we introduce discrete derivatives
of $H$, $W$, and $g$ satisfying the corresponding
exact discrete chain-rule identities. This leads
to the following lemma.
\begin{lemma}[Exact discrete chain-rule identities]
	\label{lem:discrete-chain-rules}
	For $a,b\in\mathbb R$, define
	\begin{align*}
		&\widehat{W}'(a,b)
		:=
		\frac14(a^3+a^2b+ab^2+b^3)
		-\frac12(a+b),
		\\
		&\widehat{H}'(a,b)
		:=
		\int_0^1
		\sqrt{2W\bigl(b+s(a-b)\bigr)}\,\mathrm ds,
		\\
		&\widehat{g}'(a,b)
		:=
		-\frac{\sqrt2}{2}\cos\theta
		+\frac{\sqrt2}{6}\cos\theta(a^2+ab+b^2).
	\end{align*}
	These discrete derivatives satisfy
	\begin{equation*}
		\widehat{f}'(a,b)(a-b)=f(a)-f(b),
		\qquad f=W,H,g.
	\end{equation*}
	Moreover, they are consistent with the corresponding continuous derivatives on the diagonal:
	\begin{equation*}
		\widehat{f}'(a,a)=f'(a),
		\qquad f=W,H,g.
	\end{equation*}
\end{lemma}

\begin{proof}
	For $W$ and $g$, we use the following elementary
factorizations:
	\begin{align*}
		&a^4-b^4
		=(a-b)(a^3+a^2b+ab^2+b^3),\qquad a^3-b^3
		=(a-b)(a^2+ab+b^2)\\
		&a^2-b^2
		=(a-b)(a+b).
	\end{align*}
	For $H$, the fundamental theorem of calculus gives
	\[
	H(a)-H(b)
	=(a-b)\int_0^1
	H'\bigl(b+s(a-b)\bigr)\,\mathrm ds,
	\]
	which establishes the corresponding discrete chain rule.
On the diagonal, the discrete derivatives reduce to
the corresponding ordinary derivatives:
	\[
	\widehat{W}'(a,a)=W'(a),\qquad
	\widehat{H}'(a,a)=H'(a),\qquad
	\widehat{g}'(a,a)=g'(a).
	\]
\end{proof}
The discrete derivatives of \(H\), \(W\), and \(g\) play distinct but compatible roles in the CVC construction. In particular, \(\widehat H'\) is associated with the diffuse-volume constraint, whereas \(\widehat W'\) and \(\widehat g'\) correspond to the bulk and wall energy contributions, respectively. Together with the midpoint treatment of the gradient energy and the discrete summation-by-parts identity, these constructions ensure that the fully discrete increments are compatible with the variational structure of the semi-discrete system.

In the fully discrete scheme, for $n=0,\ldots,N-1$ and $0\le j<N_x$, $0\le k\le N_y$, we use the notation
\begin{equation*}
	\widehat{W}'_{j,k}:=\widehat{W}'(U_{j,k}^{n+1},U_{j,k}^{n}),
	\qquad
	\widehat{H}'_{j,k}:=\widehat{H}'(U_{j,k}^{n+1},U_{j,k}^{n}),
	\qquad
	\widehat{g}'_j:=\widehat{g}'(U_{j,0}^{n+1},U_{j,0}^{n}).
\end{equation*}
Using these discrete derivatives, we obtain the following fully discrete scheme.
\begin{scheme}[Fully discrete CVC scheme]
	\label{sch:fully-discrete}
	Given $U^n$, find $U^{n+1}$ and $\lambda_{\eps,h}^{n+\frac12}$ such that
	\begin{equation}
		\left\{
		\begin{aligned}
			\delta_t U_{j,k}^{n+\frac12}
			&=
			\Delta_h U_{j,k}^{n+\frac12}
			-\frac{1}{\eps^2}\widehat{W}'_{j,k}
			+\lambda_{\eps,h}^{n+\frac12}\widehat{H}'_{j,k},
			\\
			\delta_t U_{j,0}^{n+\frac12}
			&=
			-\left(
			\eps(\partial_{\bn,h}U^{n+\frac12})_{j,0}
			+\widehat{g}'_j
			\right),
		\end{aligned}
		\right.
		\label{eq:fully-discrete-scheme}
	\end{equation}
	where $0\le j<N_x$, $0\le k\le N_y$, and
	\begin{equation}
		\lambda_{\eps,h}^{n+\frac12}
		=
		\frac{
			\left\langle
			\frac{1}{\eps^2}\widehat{W}'
			-\Delta_h U^{n+\frac12},
			\widehat{H}'
			\right\rangle_{\Omega,h}
		}{
			\|\widehat{H}'\|_{\Omega,h}^2
		}.
		\label{eq:fully-discrete-lambda}
	\end{equation}
\end{scheme}
\begin{remark}[Elimination of ghost values]\label{rem:ghost-elimination}
	For the implementation of the fully discrete scheme, it is convenient to eliminate the ghost values and express the evolution equations solely in terms of the physical grid values. Combining the bulk equation at $k=0$ with the dynamic boundary condition in \eqref{eq:fully-discrete-scheme} yields the following equivalent formulation:
	\begin{equation}
		\left\{
		\begin{aligned}
			\delta_t U_{j,k}^{n+\frac12}
			={}&
			\Delta_h U_{j,k}^{n+\frac12}
			-\frac{1}{\eps^2}\widehat W'_{j,k}
			+\lambda_{\eps,h}^{n+\frac12}\widehat H'_{j,k},
			&&1\le k\le N_y,\\
			\delta_t U_{j,0}^{n+\frac12}
			={}&
			\frac{1}{1+\frac{\eps\Delta y}{2}}
			\Bigg[
			\frac{\eps}{\Delta y}
			\left(U_{j,1}^{n+\frac12}
			-U_{j,0}^{n+\frac12}\right)
			\\
			{}&+
			\frac{\eps\Delta y}{2}
			\left(
			\delta_x^{\langle2\rangle}U_{j,0}^{n+\frac12}
			-\frac{1}{\eps^2}\widehat W'_{j,0}
			+\lambda_{\eps,h}^{n+\frac12}\widehat H'_{j,0}
			\right)
			-\widehat g'_j
			\Bigg],
		\end{aligned}
		\right.
		\label{eq:ghost-eliminated-scheme}
	\end{equation}
	for $0\le j<N_x$, with periodicity in the $x$-direction and the prescribed homogeneous Neumann condition at the upper boundary.
\end{remark}
\begin{remark}[Implementation of the Lagrange multiplier]
	The Lagrange multiplier can be handled in two ways. One is to eliminate the ghost values from the scheme and the Lagrange multiplier formula, and then substitute the resulting expression for $\lambda_{\eps,h}^{n+\frac12}$ into \eqref{eq:ghost-eliminated-scheme}. Alternatively, one can treat
	$\lambda_{\eps,h}^{n+\frac12}$ as an additional
	unknown and impose the discrete volume constraint
	\[
	\mathcal V_h(U^{n+1})=\mathcal V_h(U^n).
	\]
	By the exact discrete chain rule and the bulk
	equation, this formulation is equivalent to
	\eqref{eq:fully-discrete-lambda}, provided that
	$\|\widehat H'\|_{\Omega,h}>0$.
	The latter formulation is adopted in our
	numerical computations.
\end{remark}

	
	\subsection{Discrete conservation and dissipation}
	We establish the discrete volume conservation and energy dissipation properties of the fully discrete CVC scheme.
	\begin{theorem}[Discrete diffuse-volume conservation]
		\label{thm:discrete-conservation}
		The fully discrete scheme \ref{sch:fully-discrete} satisfies
		\begin{equation}
			\mathcal V_h(U^{n+1})=\mathcal V_h(U^n),
			\qquad n=0,\ldots,N-1.
			\label{eq:discrete-volume-conservation}
		\end{equation}
	\end{theorem}
	
	\begin{proof}
		The discrete chain rule in Lemma~\ref{lem:discrete-chain-rules} gives
		\begin{equation*}
			\mathcal V_h(U^{n+1})-\mathcal V_h(U^n)
			=
			\tau\left\langle
			\widehat{H}',\delta_t U^{n+\frac12}
			\right\rangle_{\Omega,h}.
		\end{equation*}
		Taking the discrete inner product of the bulk equation in \eqref{eq:fully-discrete-scheme} with $\widehat{H}'$ and using \eqref{eq:fully-discrete-lambda} gives
		\begin{equation}
			\left\langle
			\widehat{H}',\delta_t U^{n+\frac12}
			\right\rangle_{\Omega,h}=0.
			\label{eq:volume}
		\end{equation}
		Therefore, the volume conservation \eqref{eq:discrete-volume-conservation} follows immediately.
	\end{proof}
	
	\begin{theorem}[Exact discrete energy dissipation]
		\label{thm:discrete-energy-dissipation}
		The fully discrete scheme \ref{sch:fully-discrete} satisfies
		\begin{equation*}
			\mathcal E_{\eps,h}(U^{n+1})
			-\mathcal E_{\eps,h}(U^n)
			=
			-\frac{\eps}{\tau}
			\|U^{n+1}-U^n\|_{\Omega,h}^2
			-\frac{1}{\tau}
			\|U^{n+1}-U^n\|_{\Gamma,h}^2
			\le 0,
		\end{equation*}
		for $n=0,\ldots,N-1$.
	\end{theorem}
	
	\begin{proof}
		By Lemma~\ref{lem:discrete-chain-rules} and the definition of $U^{n+\frac12}$, we have
		\begin{align*}
			&\mathcal E_{\eps,h}(U^{n+1})
			-\mathcal E_{\eps,h}(U^n)
			\\
			={}&
			\eps\tau\left\langle
			\nabla_h U^{n+\frac12},
			\nabla_h\delta_t U^{n+\frac12}
			\right\rangle_h
			+\frac{\tau}{\eps}
			\left\langle
			\widehat{W}',\delta_t U^{n+\frac12}
			\right\rangle_{\Omega,h}+
			\tau\left\langle
			\widehat{g}',\delta_t U^{n+\frac12}
			\right\rangle_{\Gamma,h}.
		\end{align*}
		By the summation-by-parts identity \eqref{eq:discrete-sbp}, this becomes
		\begin{align*}
			&\mathcal E_{\eps,h}(U^{n+1})
			-\mathcal E_{\eps,h}(U^n)
			\\
			={}&
			\tau\left\langle
			-\eps\Delta_h U^{n+\frac12}
			+\frac{1}{\eps}\widehat{W}',
			\delta_t U^{n+\frac12}
			\right\rangle_{\Omega,h}+
			\tau\left\langle
			\eps\partial_{\bn,h}U^{n+\frac12}
			+\widehat{g}',
			\delta_t U^{n+\frac12}
			\right\rangle_{\Gamma,h}.
		\end{align*}
		The bulk equation in \eqref{eq:fully-discrete-scheme} gives
		\begin{align*}
			&\tau\left\langle
			-\eps\Delta_h U^{n+\frac12}
			+\frac{1}{\eps}\widehat{W}',
			\delta_t U^{n+\frac12}
			\right\rangle_{\Omega,h}\\
			={}&
			-\eps\tau
			\left\|\delta_t U^{n+\frac12}\right\|_{\Omega,h}^2
			+
			\eps\tau\lambda_{\eps,h}^{n+\frac12}
			\left\langle
			\widehat{H}',\delta_t U^{n+\frac12}
			\right\rangle_{\Omega,h}.
		\end{align*}
	The last term vanishes by \eqref{eq:volume}, so the volume-constraint correction makes no contribution to the discrete energy balance, as in the continuous system. The dynamic boundary equation gives 
		\begin{equation*}
			\eps\partial_{\bn,h}U^{n+\frac12}
			+\widehat{g}'
			=-\delta_t U^{n+\frac12}
			\qquad\text{on }\Gamma.
		\end{equation*}
		Substitution into the energy difference yields
		\begin{equation*}
			\mathcal E_{\eps,h}(U^{n+1})
			-\mathcal E_{\eps,h}(U^n)=
			-\eps\tau
			\left\|\delta_t U^{n+\frac12}\right\|_{\Omega,h}^2
			-\tau
			\left\|\delta_t U^{n+\frac12}\right\|_{\Gamma,h}^2.
		\end{equation*}
	\end{proof}

	\subsection{Energy coercivity and a mesh dependent maximum norm estimate}
	\label{sec:unconditional-energy-bounds}
	Let $U^n=(U_{j,k}^n)$ denote the approximation defined by Scheme~\ref{sch:fully-discrete}. Thus
	\begin{align}
		\left\{
		\begin{aligned}
			\frac{U^{n+1}-U^n}{\tau}
			&=\Delta_h U^{n+\frac12}
			-\frac1{\eps^2}\widehat W'
			+\lambda_{\eps,h}^{n+\frac12}\widehat H',
			&&\text{on }\Omega_h,\\
			\frac{U_{j,0}^{n+1}-U_{j,0}^n}{\tau}
			&=-\left(
			\eps(\partial_{\bn,h}U^{n+\frac12})_{j,0}
			+\widehat g'_j\right),
			&&0\leq j<N_x,
		\end{aligned}
		\right.
		\label{eq:fully-discrete-energy-scheme}
	\end{align}
	where
	\begin{align}
		U^{n+\frac12}
		=\frac{U^{n+1}+U^n}{2}, \qquad \lambda_{\eps,h}^{n+\frac12}
		=-\frac{
			\left\langle
			\Delta_hU^{n+\frac12}-\eps^{-2}\widehat W',
			\widehat H'\right\rangle_{\Omega,h}}
		{\|\widehat H'\|_{\Omega,h}^2}.
		\label{eq:fully-discrete-energy-multiplier}
	\end{align}
	Here the discrete derivatives from \cref{lem:discrete-chain-rules} are evaluated pointwise at $(U^{n+1},U^n)$, and the lower ghost values are coupled to the boundary equation. The formulas for $g$ and its divided difference are evaluated at the local wall angle. Set $c_\theta=\max_{0\le j<N_x}|\cos\theta(x_j)|$. Throughout this subsection, let
	\begin{align*}
		W(s)&=\frac14(s^2-1)^2,\qquad
		g(x,s)=\frac{\sqrt2}{6}\cos\theta(x)\,(s^3-3s),\\
		C_\theta&=\frac{\sqrt2}{3}c_\theta\,(b-a),\qquad
		\eta_\eps
		=1-c_\theta
		\left(1+\frac{\eps}{\sqrt2(d-c)}\right),\quad 
		\|U\|_{\infty,h}
		:=\max_{\substack{0\leq j<N_x\\0\leq k\leq N_y}}
		|U_{j,k}|.
	\end{align*}
	
	\begin{lemma}[Coercivity of the discrete energy]
		\label{lem:discrete-energy-coercivity}
		For every real grid function $U$,
		\begin{align}
			\eta_\eps
			\left(
			\frac{\eps}{2}\|\nabla_hU\|_h^2
			+\frac1\eps\langle W(U),1\rangle_{\Omega,h}
			\right)
			\leq \mathcal E_{\eps,h}(U)+C_\theta.
			\label{eq:discrete-energy-coercivity}
		\end{align}
		In particular, if $\eta_\eps>0$, the full discrete energy controls its nonnegative bulk part.
	\end{lemma}
	
	\begin{proof}
		Set
		\begin{align*}
			V_{j,k}=(|U_{j,k}|-1)_+,\qquad
			F(r)=\frac1{\sqrt2}\left(r^2+\frac{r^3}{3}\right),
			\qquad r\ge0.
		\end{align*}
		At each wall node, the notation $g(s)$ refers to its local wall energy. For every $s\in\mathbb R$,
		\begin{align}
			g(s)
			&\ge-\frac{\sqrt2}{3}c_\theta
			-c_\theta F((|s|-1)_+),\nonumber\\
			F((|s|-1)_+)
			&\leq\frac1{\sqrt2}W(s),\qquad
			\bigl[F'((|s|-1)_+)\bigr]^2\le2W(s).
			\label{eq:overshoot-scalar-estimates}
		\end{align}
		For $|s|\le1$, these estimates follow from $|s^3-3s|\le2$ and $F(0)=F'(0)=0$. For $s=\sigma(1+r)$, where $r\ge0$ and $\sigma\in\{-1,1\}$, one has
		\begin{gather*}
			g(s)=\sigma\cos\theta(x_j)
			\left(-\frac{\sqrt2}{3}+F(r)\right),\quad
			W(s)=\frac14r^2(r+2)^2,\quad
			[F'(r)]^2=2W(s),\\
			\frac1{\sqrt2}W(s)-F(r)
			=\frac{r^2}{\sqrt2}
			\left(\frac{r^2}{4}+\frac{2r}{3}\right)\ge0.
		\end{gather*}
		This proves \eqref{eq:overshoot-scalar-estimates} globally.
		
		Since $F'$ is nonnegative and convex on $[0,\infty)$ and $s\mapsto(|s|-1)_+$ is $1$-Lipschitz,
		\begin{align*}
			|F(V_{j,k+1})-F(V_{j,k})|
			\leq\frac12|U_{j,k+1}-U_{j,k}|
			\bigl(F'(V_{j,k+1})+F'(V_{j,k})\bigr).
		\end{align*}
		Young's inequality and \eqref{eq:overshoot-scalar-estimates} give
		\begin{align*}
			\Delta x\,|F(V_{j,k+1})-F(V_{j,k})|
			&\leq
			\frac{\eps\Delta x}{2\Delta y}
			|U_{j,k+1}-U_{j,k}|^2
			+\frac{\Delta x\Delta y}{8\eps}
			\bigl(F'(V_{j,k+1})+F'(V_{j,k})\bigr)^2
			\\
			&\quad\leq
			\frac{\eps\Delta x}{2\Delta y}
			|U_{j,k+1}-U_{j,k}|^2
			+\frac{\Delta x\Delta y}{2\eps}
			\bigl(W(U_{j,k+1})+W(U_{j,k})\bigr).
		\end{align*}
		Summing over all vertical edges, each interior potential value is counted twice and each endpoint value once. The trapezoidal weights therefore yield
		\begin{align}
			\Delta x\sum_{j=0}^{N_x-1}\sum_{k=0}^{N_y-1}
			|F(V_{j,k+1})-F(V_{j,k})|
			\leq
			\frac{\eps}{2}\|\nabla_hU\|_h^2
			+\frac1\eps\langle W(U),1\rangle_{\Omega,h}.
			\label{eq:overshoot-variation-estimate}
		\end{align}
		
		For each $j$ and $0\leq k\leq N_y$,
		\begin{align*}
			F(V_{j,0})
			\leq F(V_{j,k})
			+\sum_{\ell=0}^{N_y-1}
			|F(V_{j,\ell+1})-F(V_{j,\ell})|.
		\end{align*}
		The weights $\Delta y\,\omega_k/(d-c)$ sum to one. Averaging in $k$, summing in $j$ with weight $\Delta x$, and using \eqref{eq:overshoot-scalar-estimates} and \eqref{eq:overshoot-variation-estimate}, we obtain
		\begin{align}
			\langle F(V),1\rangle_{\Gamma,h}
			&\leq
			\frac1{d-c}\langle F(V),1\rangle_{\Omega,h}+
			\Delta x\sum_{j=0}^{N_x-1}\sum_{k=0}^{N_y-1}
			|F(V_{j,k+1})-F(V_{j,k})|
			\nonumber\\
			&\leq
			\frac{\eps}{2}\|\nabla_hU\|_h^2
			+\left(\frac1\eps+\frac1{\sqrt2(d-c)}\right)
			\langle W(U),1\rangle_{\Omega,h}
			\nonumber\\
			&\leq
			\left(1+\frac{\eps}{\sqrt2(d-c)}\right)
			\left(
			\frac{\eps}{2}\|\nabla_hU\|_h^2
			+\frac1\eps\langle W(U),1\rangle_{\Omega,h}
			\right).
			\label{eq:overshoot-trace-estimate}
		\end{align}
		Finally, the definition \eqref{eq:discrete-energy} and the first inequality in \eqref{eq:overshoot-scalar-estimates} imply
		\begin{align*}
			\mathcal E_{\eps,h}(U)+C_\theta
			&\ge
			\frac{\eps}{2}\|\nabla_hU\|_h^2
			+\frac1\eps\langle W(U),1\rangle_{\Omega,h}
			-c_\theta\langle F(V),1\rangle_{\Gamma,h}.
		\end{align*}
		Substituting \eqref{eq:overshoot-trace-estimate} proves \eqref{eq:discrete-energy-coercivity}.
	\end{proof}
	
	\begin{theorem}
		\label{thm:fully-discrete-energy-bounds}
		Let $\tau>0$ and let $(U^n)_{n=0}^N$ solve \eqref{eq:fully-discrete-energy-scheme}--\eqref{eq:fully-discrete-energy-multiplier}, with $\|\widehat H'\|_{\Omega,h}>0$ at every step. Then $\mathcal V_h(U^n)=\mathcal V_h(U^0)$ and
		\begin{align}
			\mathcal E_{\eps,h}(U^n)
			&+\sum_{\ell=0}^{n-1}
			\left(
			\frac{\eps}{\tau}
			\|U^{\ell+1}-U^\ell\|_{\Omega,h}^2
			+\frac1\tau
			\|U^{\ell+1}-U^\ell\|_{\Gamma,h}^2
			\right)
			=\mathcal E_{\eps,h}(U^0).
			\label{eq:fully-discrete-energy-identity}
		\end{align}
		If, in addition, $\eta_\eps>0$, then
		\begin{align}
			&\eta_\eps
			\left(
			\frac{\eps}{2}\|\nabla_hU^n\|_h^2
			+\frac1\eps\langle W(U^n),1\rangle_{\Omega,h}
			\right)
			\nonumber\\
			&\quad+
			\sum_{\ell=0}^{n-1}
			\left(
			\frac{\eps}{\tau}
			\|U^{\ell+1}-U^\ell\|_{\Omega,h}^2
			+\frac1\tau
			\|U^{\ell+1}-U^\ell\|_{\Gamma,h}^2
			\right)
			\leq \mathcal E_{\eps,h}(U^0)+C_\theta,
			\label{eq:fully-discrete-energy-bound}
		\end{align}
		and
		\begin{align}
			\|U^n\|_{\infty,h}
			\le
			\left[
			1+2\sqrt{
				\frac{2\eps\bigl(\mathcal E_{\eps,h}(U^0)+C_\theta\bigr)}
				{\eta_\eps\Delta x\Delta y}}
			\right]^{1/2},
			\qquad 0\leq n\leq N.
			\label{eq:fully-discrete-maximum-norm-bound}
		\end{align}
		These estimates impose no restriction on $\tau$ and require no initial pointwise bound.
	\end{theorem}
	
	\begin{proof}
		Taking the bulk inner product in \eqref{eq:fully-discrete-energy-scheme} with $\widehat H'$ and using \eqref{eq:fully-discrete-energy-multiplier}, we obtain
		\begin{align*}
			\langle\widehat H',U^{n+1}-U^n\rangle_{\Omega,h}=0.
		\end{align*}
		The discrete chain rule therefore gives $\mathcal V_h(U^{n+1})=\mathcal V_h(U^n)$. Testing the bulk equation by $\eps(U^{n+1}-U^n)$, applying \eqref{eq:discrete-sbp}, and using the boundary equation gives
		\begin{align*}
			\frac{\eps}{\tau}
			\|U^{n+1}-U^n\|_{\Omega,h}^2
			&=-\eps
			\langle\nabla_hU^{n+\frac12},
			\nabla_h(U^{n+1}-U^n)\rangle_h
			-\frac1\eps
			\langle\widehat W',U^{n+1}-U^n\rangle_{\Omega,h}\\
			&\quad-\frac1\tau
			\|U^{n+1}-U^n\|_{\Gamma,h}^2
			-\langle\widehat g',U^{n+1}-U^n\rangle_{\Gamma,h}.
		\end{align*}
		The discrete chain rules and
		\begin{align*}
			\langle\nabla_hU^{n+\frac12},
			\nabla_h(U^{n+1}-U^n)\rangle_h
			=\frac12\left(
			\|\nabla_hU^{n+1}\|_h^2-\|\nabla_hU^n\|_h^2
			\right)
		\end{align*}
		therefore imply
		\begin{align*}
			\mathcal E_{\eps,h}(U^{n+1})-\mathcal E_{\eps,h}(U^n)
			=-\frac{\eps}{\tau}
			\|U^{n+1}-U^n\|_{\Omega,h}^2
			-\frac1\tau
			\|U^{n+1}-U^n\|_{\Gamma,h}^2.
		\end{align*}
		Summation proves \eqref{eq:fully-discrete-energy-identity}. Combining this identity with \cref{lem:discrete-energy-coercivity} gives \eqref{eq:fully-discrete-energy-bound}.
		
		For each physical node, $\omega_k\ge1/2$, and hence
		\begin{align*}
			\frac{\Delta x\Delta y}{8\eps}
			\bigl((U_{j,k}^n)^2-1\bigr)^2
			\leq\frac{\Delta x\Delta y\,\omega_k}{\eps}
			W(U_{j,k}^n)
			\leq\frac1\eps\langle W(U^n),1\rangle_{\Omega,h}
			\leq
			\frac{\mathcal E_{\eps,h}(U^0)+C_\theta}{\eta_\eps}.
		\end{align*}
		Consequently,
		\begin{align*}
			|U_{j,k}^n|^2
			\le
			1+2\sqrt{
				\frac{2\eps\bigl(\mathcal E_{\eps,h}(U^0)+C_\theta\bigr)}
				{\eta_\eps\Delta x\Delta y}}.
		\end{align*}
		Taking square roots and maximizing over the grid proves \eqref{eq:fully-discrete-maximum-norm-bound}.
	\end{proof}
	
	\begin{remark}
		\label{rem:energy-bound-mesh-dependence}
		The energy identity \eqref{eq:fully-discrete-energy-identity} holds for every $\tau>0$. The coercive estimate additionally requires $\eta_\eps>0$, which holds for $c_\theta<1$ when $\eps$ is sufficiently small. If the wall is neutral, then $c_\theta=0$ and $\eta_\eps=1$ for every $\eps>0$. The same estimates apply to spatially varying wall angles with constants depending on the maximum of $|\cos\theta(x_j)|$.
		
		Suppose that $\mathcal E_{\eps,h}(U^0)+C_\theta$ is bounded and $\eta_\eps$ is bounded away from zero. Then the right-hand side of \eqref{eq:fully-discrete-maximum-norm-bound} tends to $1$ if $\eps/(\Delta x\Delta y)\to0$. For fixed mesh sizes this gives an excess of order $\sqrt\eps$. For $\Delta x=\Delta y=\eps$, the estimate does not yield a bound tending to $1$. Thus it is an unconditional energy-based maximum norm estimate with explicit mesh dependence.
	\end{remark}
	
	\section{Numerical experiments}\label{sec:numerical}
	
	\subsection{Accuracy verification}
	We test the CVC scheme on $\Omega=[0,2\pi]^2$ with periodic conditions in $x$, dynamic relaxation on the lower wall, and homogeneous Neumann conditions above. We use $\eps=1$, $\theta_Y=2\pi/3$, $u_0=\cos x\cos y$, and $T=1$. For temporal refinement, $N_x=N_y=128$ and $\tau=1/20,1/40,1/80,1/160,1/320$; errors compare successive time steps. For spatial refinement, $\tau=10^{-3}$ and $N_x=N_y=N=8,16,32,64,128$; errors compare successive grids after restriction to coarse nodes. Consecutive errors give the orders in Tables~\ref{tab:convergence-order}, confirming second-order convergence in both $L^2$ and $L^\infty$.
	
	\begin{table}[!t]
	\centering
	\caption{Temporal and spatial convergence results of the proposed scheme.}
	\label{tab:convergence-order}
	\setlength{\tabcolsep}{4pt}
	\resizebox{\linewidth}{!}{%
	\begin{tabular}{ccccc@{\hspace{12pt}}ccccc}
		\toprule
		\multicolumn{5}{c}{Temporal convergence} & \multicolumn{5}{c}{Spatial convergence} \\
		\cmidrule(r{12pt}){1-5} \cmidrule(l){6-10}
		$\tau$ & $L^2$ error & Order & $L^\infty$ error & Order & $\Delta x$ & $L^2$ error & Order & $L^\infty$ error & Order \\
		\midrule
		$1/20$  & $3.6106\times10^{-4}$ & --       & $1.3040\times10^{-4}$ & --       & $\pi/4$  & $5.9233\times10^{-2}$ & --       & $1.9368\times10^{-2}$ & --       \\
		$1/40$  & $8.9573\times10^{-5}$ & $2.0111$ & $3.2584\times10^{-5}$ & $2.0007$ & $\pi/8$  & $1.4289\times10^{-2}$ & $2.0515$ & $5.0171\times10^{-3}$ & $1.9487$ \\
		$1/80$  & $2.2391\times10^{-5}$ & $2.0001$ & $8.1452\times10^{-6}$ & $2.0002$ & $\pi/16$ & $3.5367\times10^{-3}$ & $2.0144$ & $1.2611\times10^{-3}$ & $1.9921$ \\
		$1/160$ & $5.5976\times10^{-6}$ & $2.0000$ & $2.0362\times10^{-6}$ & $2.0000$ & $\pi/32$ & $8.8197\times10^{-4}$ & $2.0036$ & $3.1568\times10^{-4}$ & $1.9982$ \\
		\bottomrule
	\end{tabular}}
\end{table}
	\subsection{Improved volume preservation over the mass-conserving model}
	On the homogeneous lower wall of $\Omega=(0,1)\times(0,0.5)$, we compare small-droplet retention under nonlinear-volume and linear-mass constraints, then verify cap shapes and contact angles. In two dimensions, volume is represented by area.
	
	For a prescribed sharp-interface reference area $A_{\mathrm{target}}$ and a Young angle $\theta_Y\in(0,\pi)$, the corresponding ideal circular cap has radius $R=\left(\frac{A_0}{\theta_Y-\sin\theta_Y\cos\theta_Y}\right)^{1/2}$, footprint width $\ell=2R\sin\theta_Y$, and height $h_{\mathrm{cap}}=R(1-\cos\theta_Y)$.
	This reference applies when the remaining boundaries do not constrain the cap and gravity is absent, as in the present energy.
	
	We first compare the proposed nonlinear-volume-constrained model with the mass-conserving Allen--Cahn model. The latter preserves the linear integral $M(u)=\int_\Omega u\,\mathrm dx$. We prescribe the wetting angle $\theta_Y=60^\circ$ and use initial droplets with radii $R_0=0.30$, $0.16$, and $0.05$.
	
	The computations use $\eps=0.01$, $\tau=0.05$, and a spatial grid of $128\times64$. The initial contact angle is $\theta_0=90^\circ$. The initial phase field is constructed as
		\[
		u_0(x,y)=\tanh\left(-\frac{\phi_0(x,y)-s_0}{\sqrt{2}\eps}\right),
		\qquad
		\phi_0(x,y)=\sqrt{(x-\tfrac12)^2+y^2}-R_0,
		\]
		where the scalar shift $s_0$ is chosen to match the prescribed discrete nonlinear volume associated with the target cap area $A_{\mathrm{target}}=\pi R_0^2/2$. The simulations are run up to $T_{\max}=200$, with each run allowed to terminate earlier when the near-stationarity criterion is satisfied.
	
\begin{figure}[!t]
		\centering
		\includegraphics[width=0.84\textwidth]{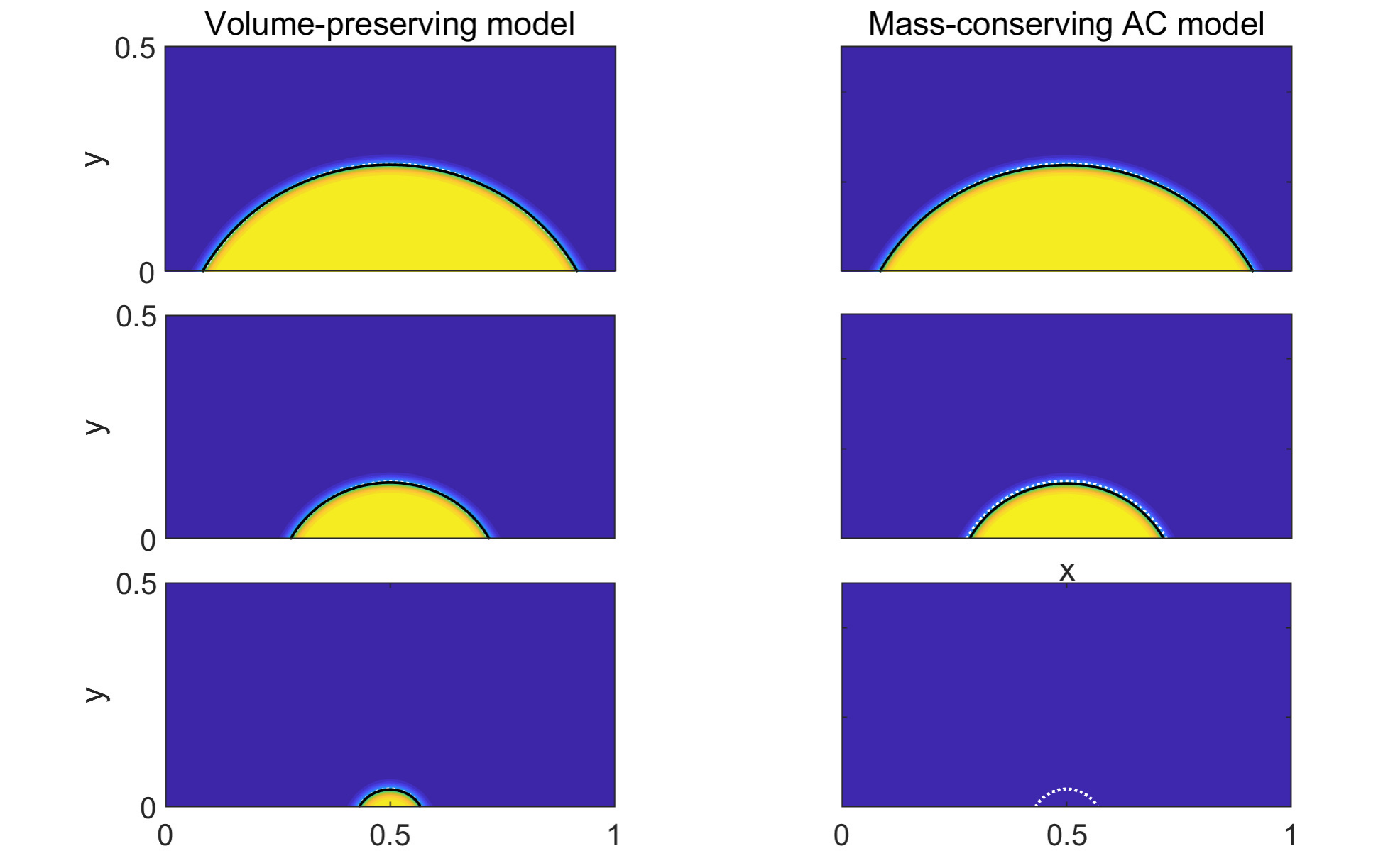}
		\caption{Proposed volume-preserving model (left) and mass-conserving Allen--Cahn model (right) on a flat wall, $\theta_Y=60^\circ$. Initial radii: $R_0=0.30$, $0.16$, $0.05$ (top to bottom); dashed curves: ideal caps of the prescribed area and angle.}
		\label{fig:reference-model-comparison}
	\end{figure}

	Figure~\ref{fig:reference-model-comparison} compares the computed shapes with ideal caps of the prescribed area and Young angle. Both models agree for the largest droplet, but the mass-conserving model increasingly underestimates size as $R_0$ decreases and loses the $R_0=0.05$ droplet. The proposed model retains it close to the reference cap, demonstrating the benefit of the nonlinear constraint for small droplets.
	
\begin{figure}[t!]
		\centering
		\includegraphics[width=0.90\textwidth]{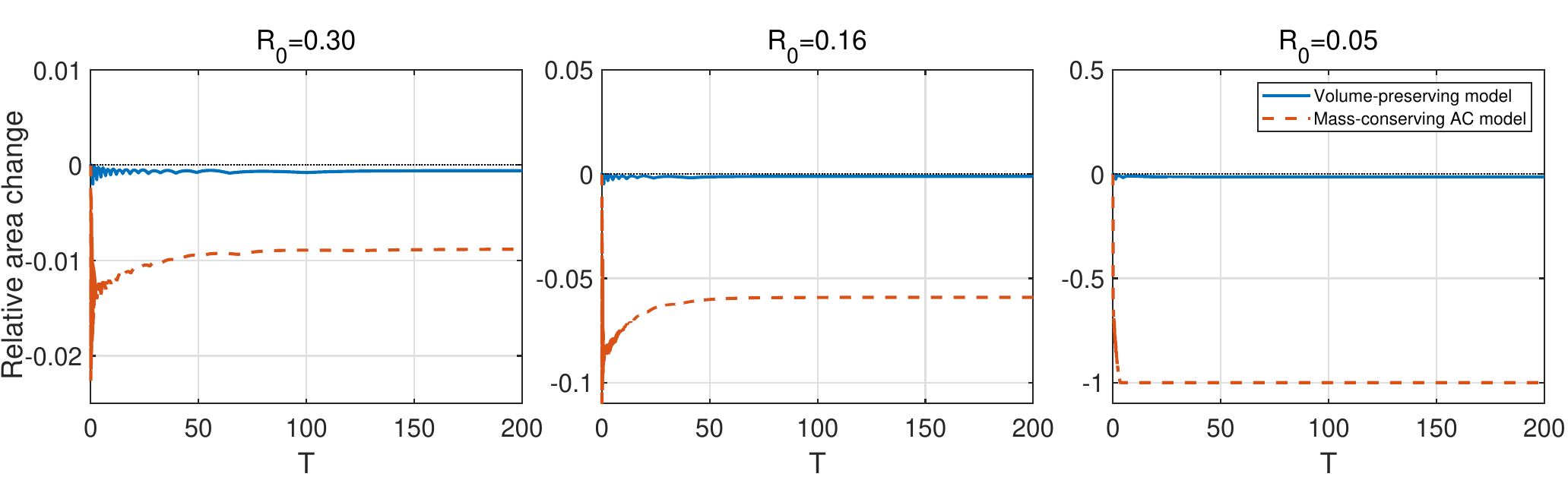}
		\caption{Relative geometric-area change $(A(t)-A_0)/A_0$ at $\theta_Y=60^\circ$ for $R_0=0.30$, $0.16$, $0.05$ (left to right). Solid: proposed model; dashed: mass-conserving model. Panel scales differ.}
		\label{fig:reference-area-change}
	\end{figure}

	Figure~\ref{fig:reference-area-change} reports $(A(t)-A_0)/A_0$. We estimate $A(t)$ by piecewise-linear interpolation of the $u=0$ interface on vertical grid edges, summing the $u>0$ fractions with weight $\Delta x\Delta y$. This geometric diagnostic differs from the conserved diffuse-volume. Its relative change remains near zero for the proposed model but reaches $-1$ when the mass-conserving model loses the smallest droplet; the largest droplet is much less affected, consistent with finite-thickness sensitivity in \cite{YueZhouFeng2007}.
	
\begin{figure}[!t]
		\centering
		\includegraphics[width=0.4\textwidth]{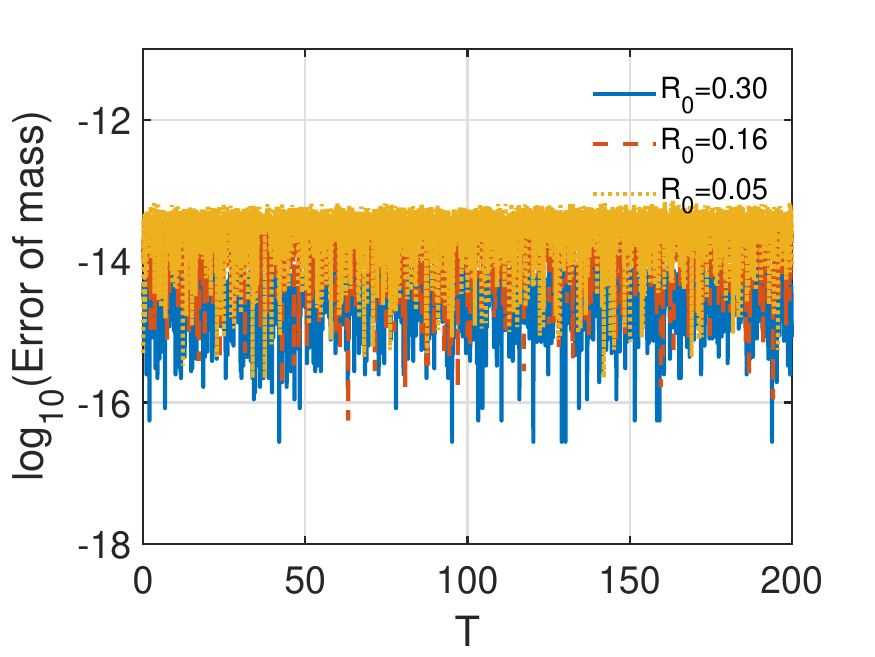}
		\caption{Absolute mass error $\log_{10}|M_h(t)-M_h(0)|$ of the mass-conserving model, $\theta_Y=60^\circ$, $R_0=0.30$, $0.16$, $0.05$. Mass is conserved even when the smallest droplet disappears (Figure~\ref{fig:reference-area-change}).}
		\label{fig:reference-mass-error}
	\end{figure}

	Figure~\ref{fig:reference-mass-error} plots $\log_{10}|M_h(t)-M_h(0)|$, the unnormalized absolute discrete error in $M(u)=\int_\Omega u\,\mathrm dx$. It stays near $10^{-13}$ or below, even when the droplet disappears. At finite interface thickness, small bulk-phase shifts over the larger surrounding region can compensate for loss of the droplet while reducing interfacial energy. This mechanism was analyzed for Cahn--Hilliard dynamics in \cite{YueZhouFeng2007}; the uniform multiplier in the mass-conserving Allen--Cahn model likewise permits bulk redistribution. In contrast, the proposed correction is weighted by $H'(u)$, which vanishes at $u=\pm1$. The comparisons demonstrate improved geometric retention; exact conservation concerns the nonlinear diffuse-volume.
	
	We next use $\theta_Y=60^\circ$ and $120^\circ$, with the same $\eps=0.01$, $\tau=0.05$, $128\times64$ grid, shifted-tanh initialization, and $\theta_0=90^\circ$. Here $R_0=0.16$, $A_{\mathrm{target}}=\pi(0.16)^2/2$, and $T_{\max}=100$, with the same near-stationarity stopping criterion. The shift $s_0$ matches the prescribed discrete nonlinear volume. We average the left and right contact angles, evaluating $\operatorname{atan2}(|u_x|,-u_y)$ at each $u=0$ contact point using second-order centered differences in $x$, second-order one-sided differences in $y$, and linear interpolation to the contact point.
	
\begin{figure}[!t]
		\centering
		\includegraphics[width=0.8\textwidth]{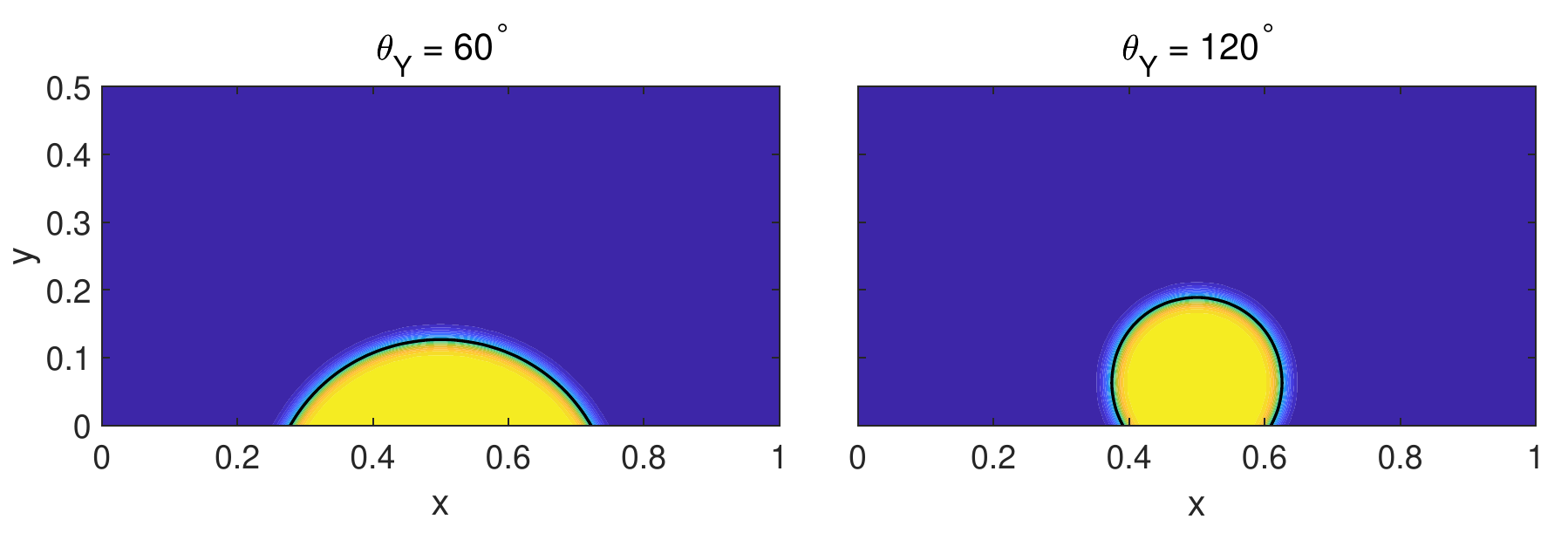}
		\caption{Flat-wall cap shapes for $\theta_Y=60^\circ$ (left) and $120^\circ$ (right); see Table~\ref{tab:contact-angle} for numerical angles and errors.}
		\label{fig:reference-contact-shapes}
	\end{figure}

	Figure~\ref{fig:reference-contact-shapes} shows the cap shapes. Table~\ref{tab:contact-angle} tests five Young angles from $30^\circ$ to $150^\circ$ on a $256\times128$ grid. The errors $100|\theta_h-\theta_Y|/\theta_Y$ range from $0.0101\%$ to $6.1252\%$, with the smallest at $90^\circ$, supporting accurate flat-wall wetting across these preferences.
	
	\begin{table}[!t]
	\centering
	\caption{Numerical contact angles and relative errors for the homogeneous flat-wall test on a $256\times128$ grid.}
	\label{tab:contact-angle}
	\begin{tabular}{lccccc}
		\toprule
		$\theta_Y$      & $30^\circ$      & $60^\circ$      & $90^\circ$      & $120^\circ$      & $150^\circ$      \\
		\midrule
		$\theta_h$      & $28.1625^\circ$ & $58.7536^\circ$ & $89.9909^\circ$ & $119.7431^\circ$ & $147.0447^\circ$ \\
		Rel. error (\%) & $6.1252$        & $2.0773$        & $0.0101$        & $0.2141$         & $1.9702$         \\
		\bottomrule
	\end{tabular}
\end{table}
	
\begin{figure}[htbp]
	\centering
	\begin{minipage}[t]{0.35\textwidth}
		\centering
		\includegraphics[width=\linewidth]{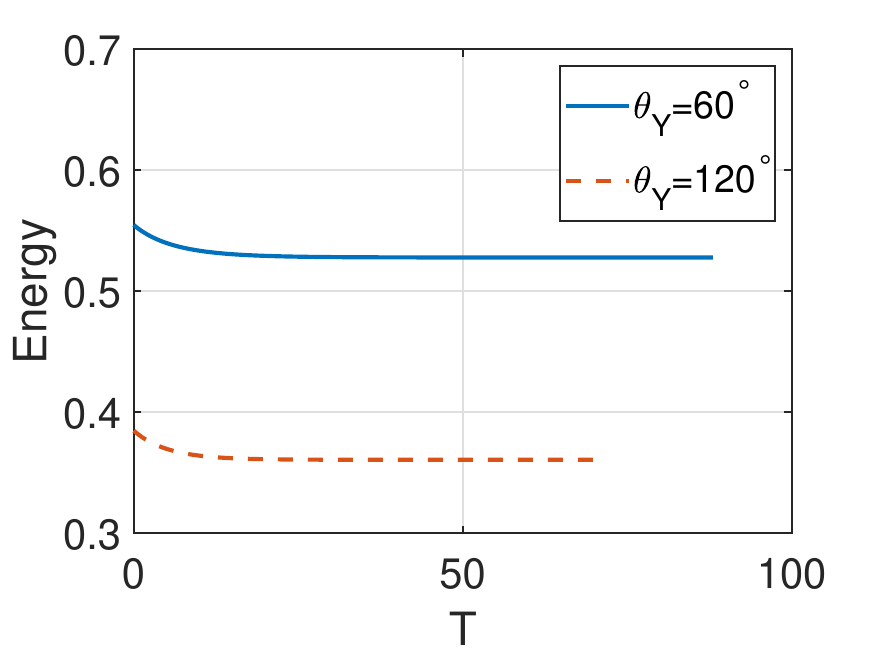}
	\end{minipage}%
	\hspace{0.02\textwidth}
	\begin{minipage}[t]{0.35\textwidth}
		\centering
		\includegraphics[width=\linewidth]{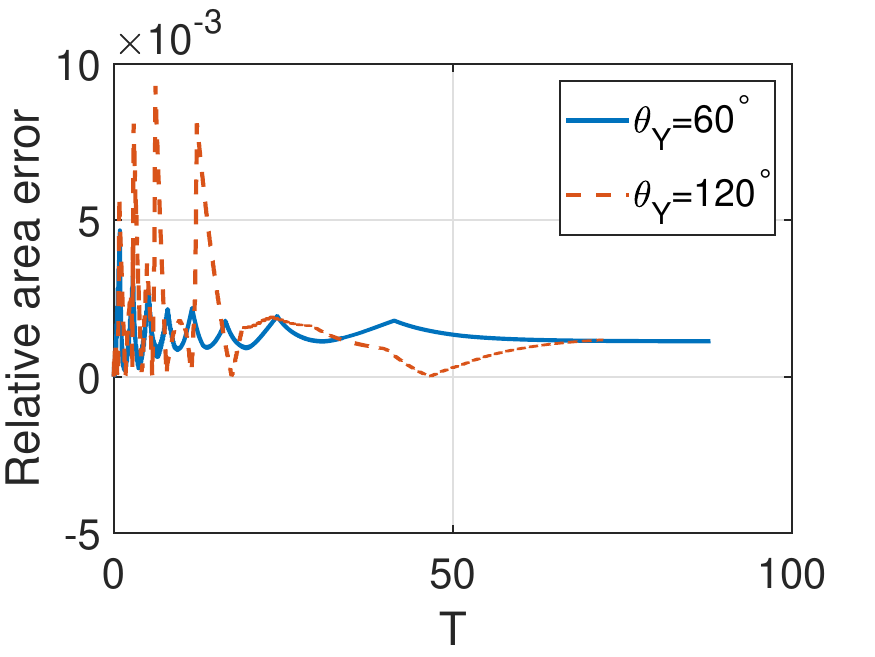}
	\end{minipage}
	\caption{Discrete energy (left) and relative geometric-area error (right) for the proposed model with $\theta_Y=60^\circ$ and $120^\circ$ on a homogeneous flat wall.}
	\label{fig:reference-energy-area}
\end{figure}

	Figure~\ref{fig:reference-energy-area} shows decreasing energy and relative geometric-area errors on the $10^{-3}$ scale for both Young angles. Together with droplet retention and contact-angle agreement, these results connect the nonlinear constraint and its compatible discretization to improved wetting computations; geometric area is not exactly conserved.
	
	\subsection{Wetting and pinning on chemically heterogeneous substrates}
	
	We test contact-line pinning through two control pairs: Cases~I/II and III/IV start with different shapes and footprints but are designed to approach the same pinning positions within each pair.
	
	The computational domain is $\Omega=[-1,1]\times[0,1]$, with a uniform grid of $N_x=256$ and $N_y=128$, giving $\Delta x=\Delta y=1/128$. We take $\eps=0.02$, the time step $\tau=10^{-2}$, and the final time $T=20$. Periodic boundary conditions are imposed in the $x$-direction, together with a dynamic boundary condition on the lower wall and a homogeneous Neumann condition on the upper boundary.
	
	The substrate consists of two materials with Young angles $\theta_A=120^\circ$ and $\theta_B=60^\circ$. Material B occupies
	$[-0.58,-0.46]\cup[-0.34,-0.22]
	\cup[0.22,0.34]\cup[0.46,0.58]$,
	while the remaining portions belong to material A. The four initial droplets are circular caps centered at $x=0$, with half-widths $d_0=0.28$, $0.42$, $0.52$, and $0.70$ for Cases~I--IV, respectively. Their initial contact angles are chosen so that all four caps have the same sharp-interface area $A_{\mathrm{target}}=0.30$. A scalar level-set adjustment is then applied to match their initial discrete diffuse-volumes. Cases~I and~II are designed to approach pinning positions near $x=\pm0.34$, whereas Cases~III and~IV target positions near $x=\pm0.58$.
	
	\begin{figure}[!t]
		\centering
		\includegraphics[width=0.90\textwidth]
		{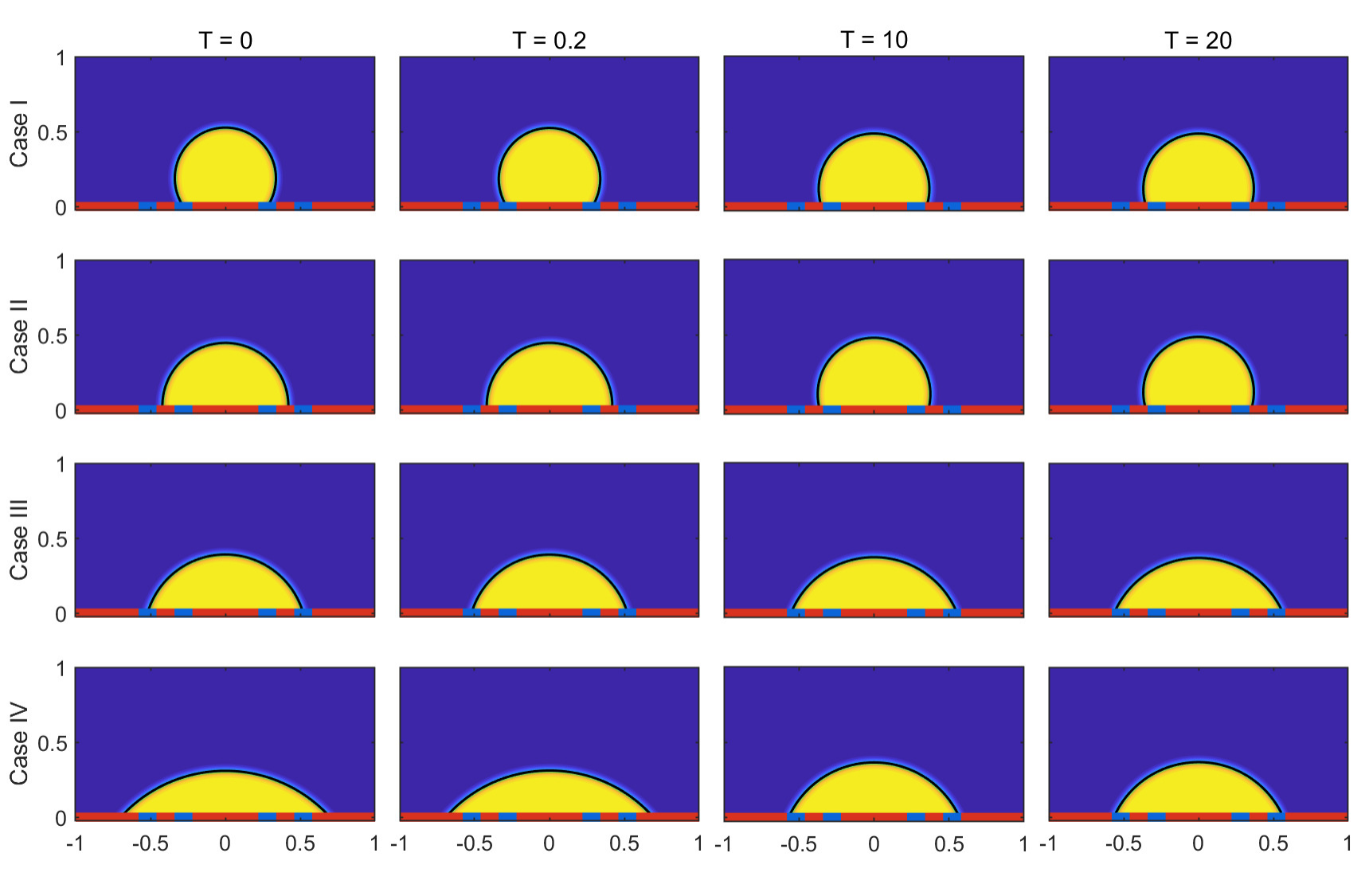} \caption{Paired pinning tests on a chemically patterned wall. Rows: Cases~I--IV; columns: $T=0$, $0.2$, $10$, $20$. Cases~I/II and III/IV approach common pinning positions within each pair. Colored wall segments mark the chemical pattern.}
		\label{fig:heterogeneous-evolution}
	\end{figure}
	
\begin{figure}[!t]
		\centering
		\begin{minipage}[t]{0.99\textwidth}
			\centering
			\includegraphics[width=0.55\textwidth]{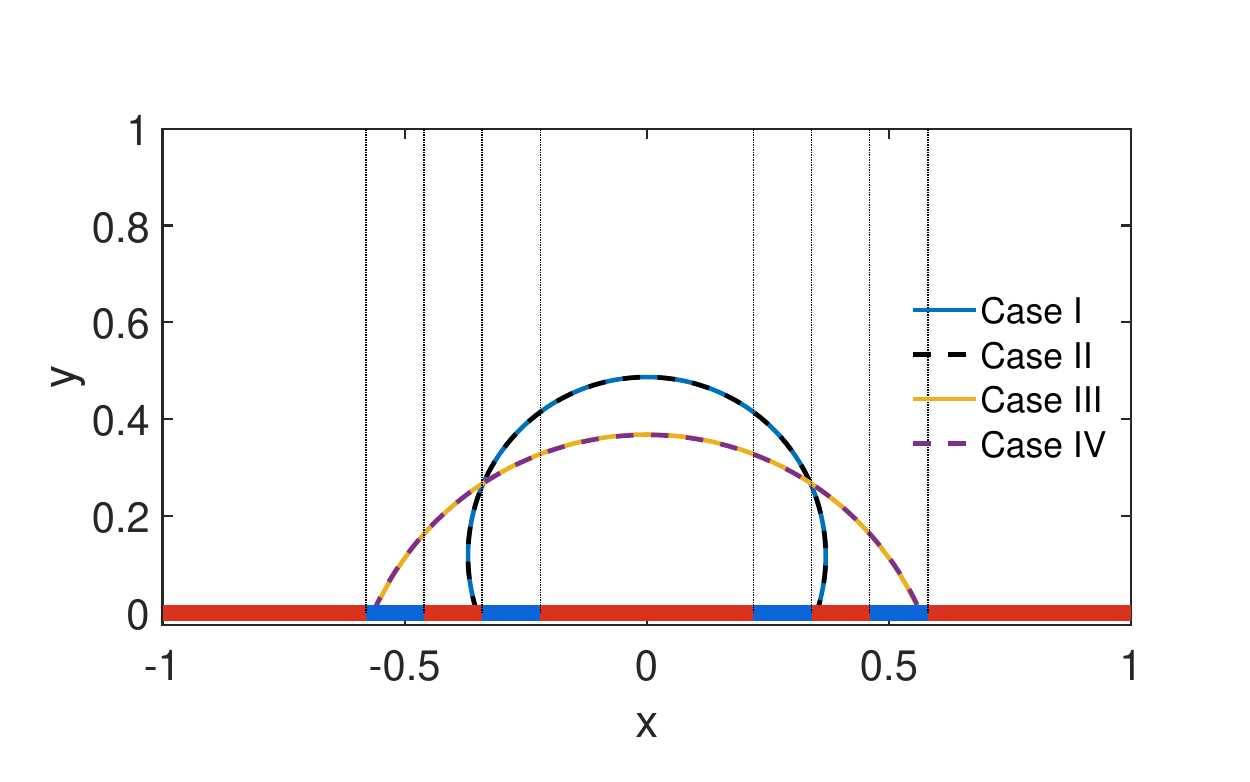}
		\end{minipage}
		\caption{Final pinned contours: Cases~I and III (solid), II and IV (dashed). Profiles nearly coincide within each pair; the two pairs select distinct configurations.}
		\label{fig:heterogeneous-final}
	\end{figure}

\begin{figure}[t!]
		\centering
		\begin{minipage}[t]{0.4\textwidth}
			\centering
			\includegraphics[width=0.8\linewidth]{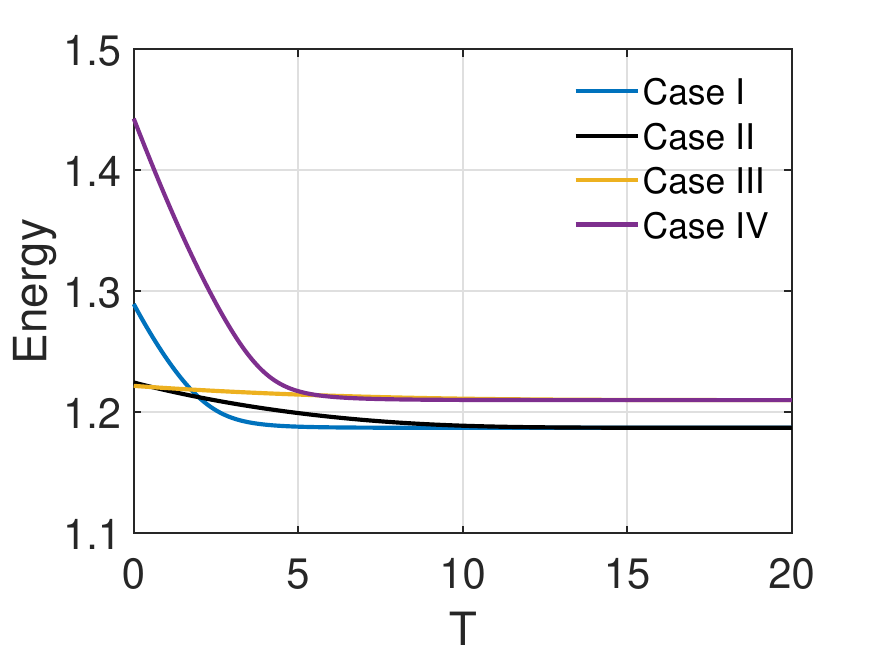}
		\end{minipage}
		\begin{minipage}[t]{0.4\textwidth}
			\centering
			\includegraphics[width=0.8\linewidth]{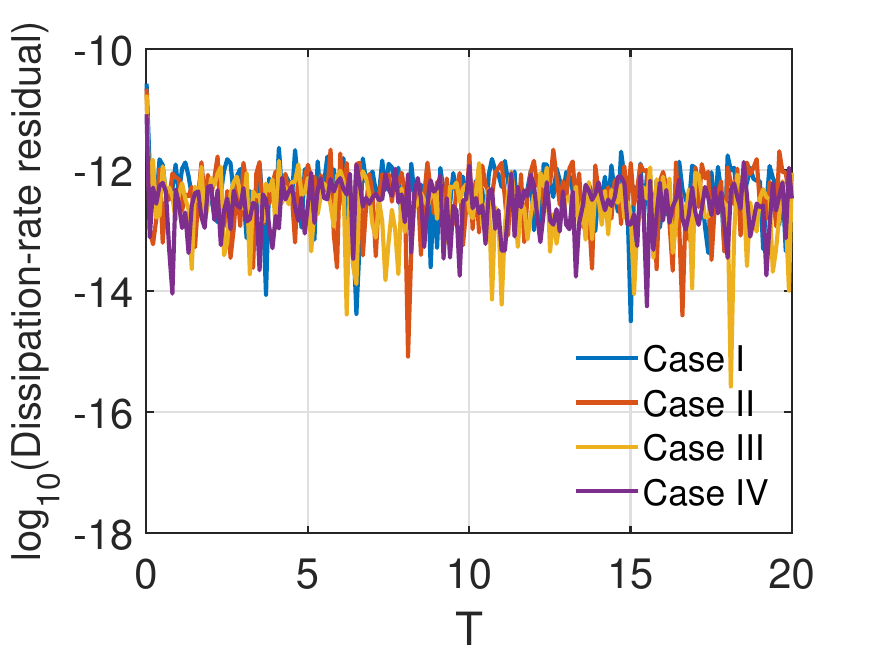}
		\end{minipage}
		\begin{minipage}[t]{0.4\textwidth}
			\centering
			\includegraphics[width=0.8\linewidth]{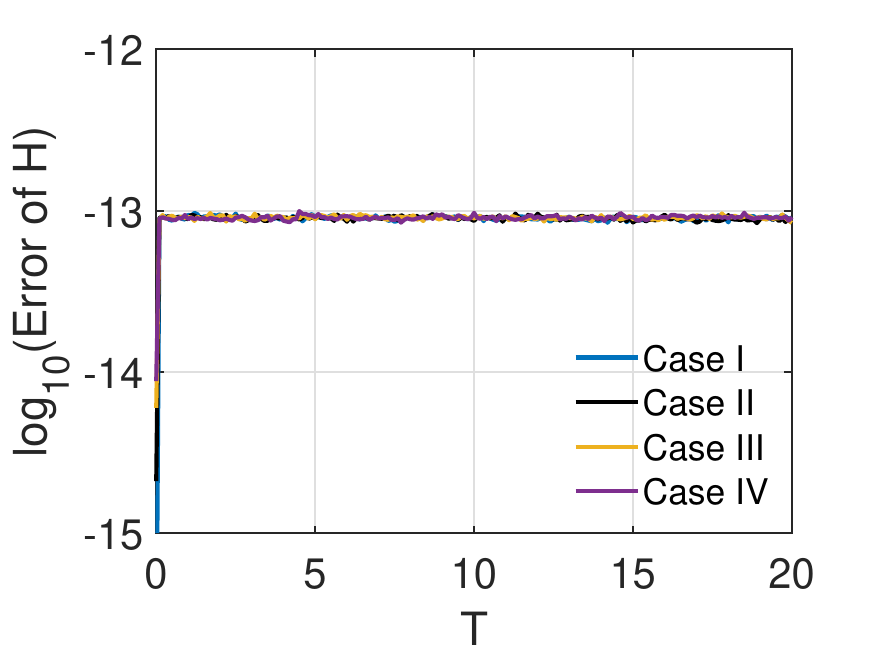}
		\end{minipage}
		\begin{minipage}[t]{0.4\textwidth}
			\centering
			\includegraphics[width=0.8\linewidth]{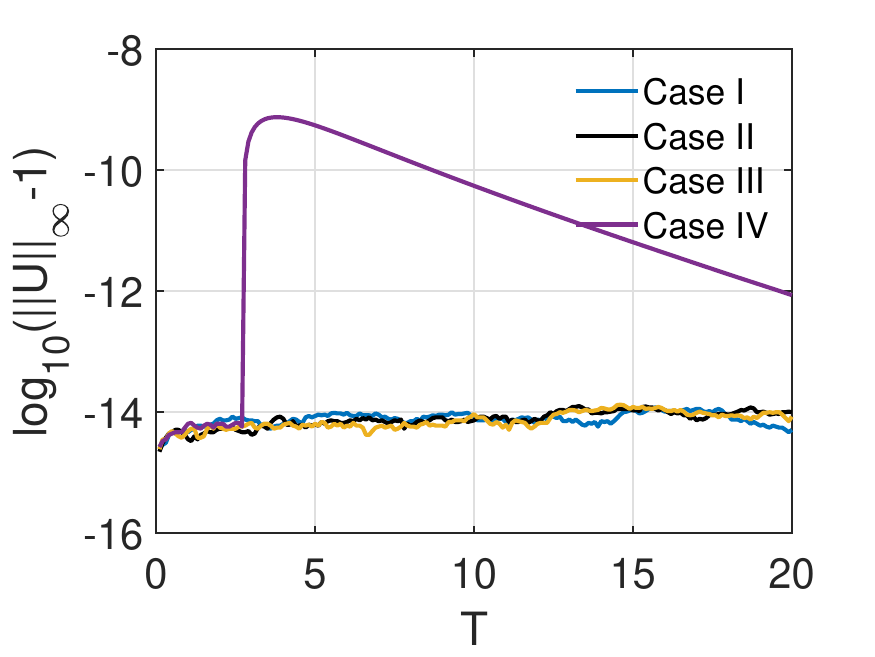}
		\end{minipage}
		\caption{Discrete energy evolution (top left), logarithmic energy-dissipation residual (top right), diffuse-volume conservation error (bottom left), and logarithmic maximum-norm diagnostic $\log_{10}|\|U\|_\infty-1|$ (bottom right) for the heterogeneous-wettability test.}
		\label{fig:heterogeneous-structure}
	\end{figure}

	Figures~\ref{fig:heterogeneous-evolution} and~\ref{fig:heterogeneous-final} show nearly coincident final profiles within each pair and distinct pinned configurations between pairs. Pinning is reproducible across the initial shapes in each pair. The diagnostics in Figure~\ref{fig:heterogeneous-structure} show energy decay, small diffuse-volume error, and maximum-norm deviations from $[-1,1]$.
	
	The scheme also preserves the original energy-dissipation rate identity:
	\begin{equation*}
		\frac{\mathcal E_{\eps,h}(U^{n+1})
			-\mathcal E_{\eps,h}(U^n)}{\tau}
		=-\mathcal D_h^{n+\frac12},
	\end{equation*}
	where
	$\mathcal D_h^{n+\frac12}
		=
		\eps\|\delta_t U^{n+\frac12}\|_{\Omega,h}^2
		+
		\|\delta_t U^{n+\frac12}\|_{\Gamma,h}^2,
		~
		\delta_t U^{n+\frac12}
		=\frac{U^{n+1}-U^n}{\tau}.
	$
	To verify this identity numerically, we define the rate residual
	\begin{equation*}
		\mathcal R_h^{n+\frac12}
		=
		\frac{\mathcal E_{\eps,h}(U^{n+1})
			-\mathcal E_{\eps,h}(U^n)}{\tau}
		+\mathcal D_h^{n+\frac12}.
	\end{equation*}
	The upper-right panel of Figure~\ref{fig:heterogeneous-structure} plots its logarithmic absolute value, which remains small to nonlinear-solver accuracy. The residual remains small up to the nonlinear solver tolerance, confirming the discrete energy-dissipation identity during the pinning process.
	
	
	\subsection{Chemically directed migration and coalescence}
	
\begin{figure}[t!]
		\centering
		\includegraphics[width=0.8\textwidth]{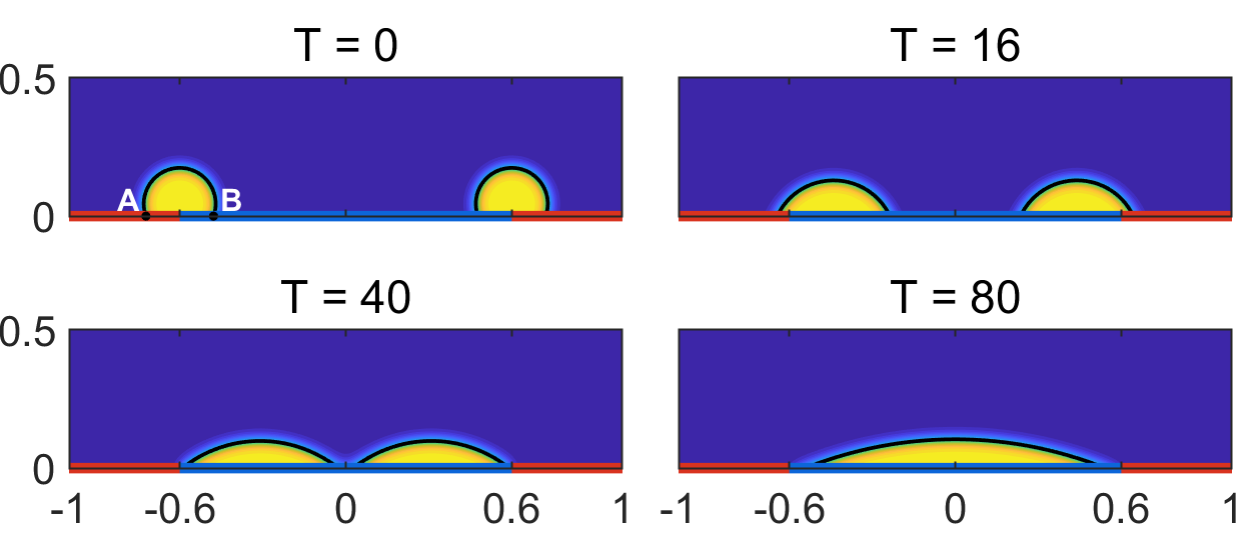}
		\caption{Evolution of two symmetric droplets on a surface with heterogeneous wettability.}
		\label{selfclean}
	\end{figure}

\begin{figure}[t!]
		\centering
		\includegraphics[width=0.7\textwidth]{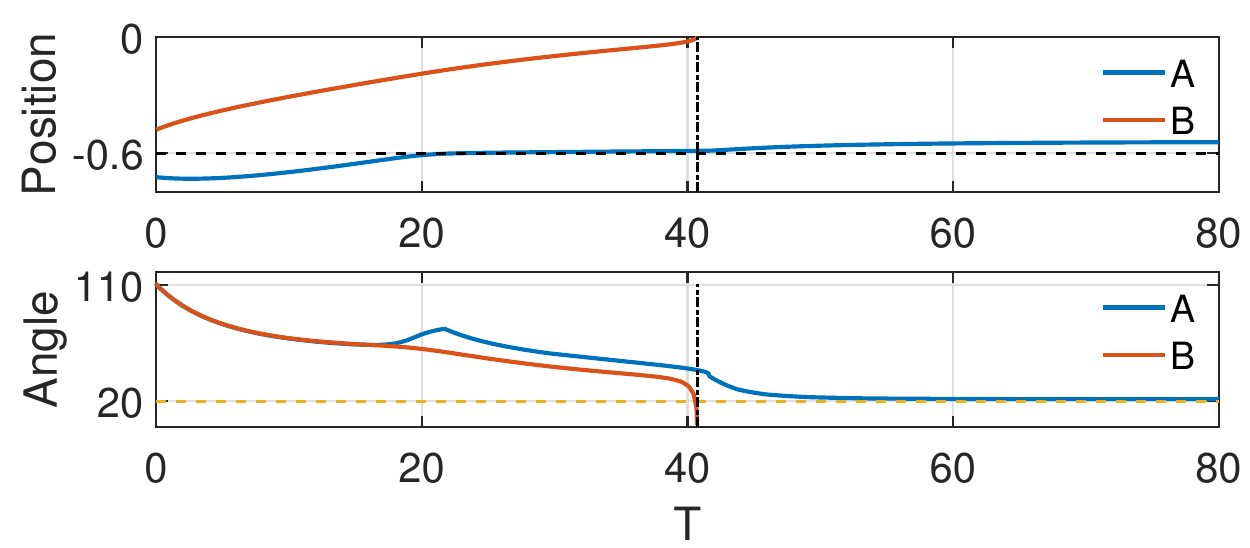}
		\caption{Contact-line positions (top) and contact angles (bottom) of the droplet in Figure~\ref{selfclean}.}
		\label{fig:migration-contact}
	\end{figure}

\begin{figure}[htbp]
		\centering
		\begin{minipage}[t]{0.32\textwidth}
			\centering
			\includegraphics[width=0.96\linewidth]{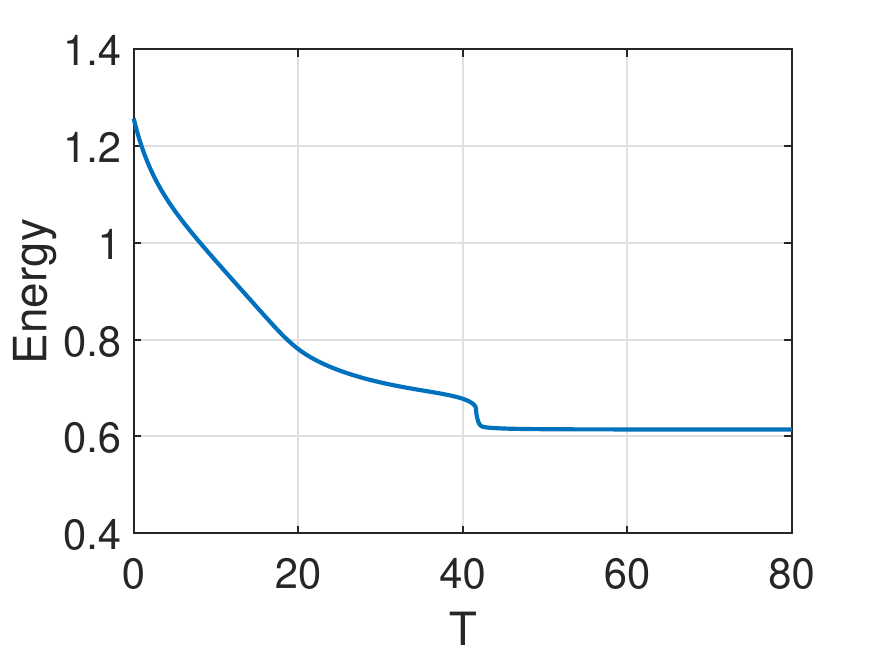}
		\end{minipage}
		\begin{minipage}[t]{0.32\textwidth}
			\centering
			\includegraphics[width=0.96\linewidth]{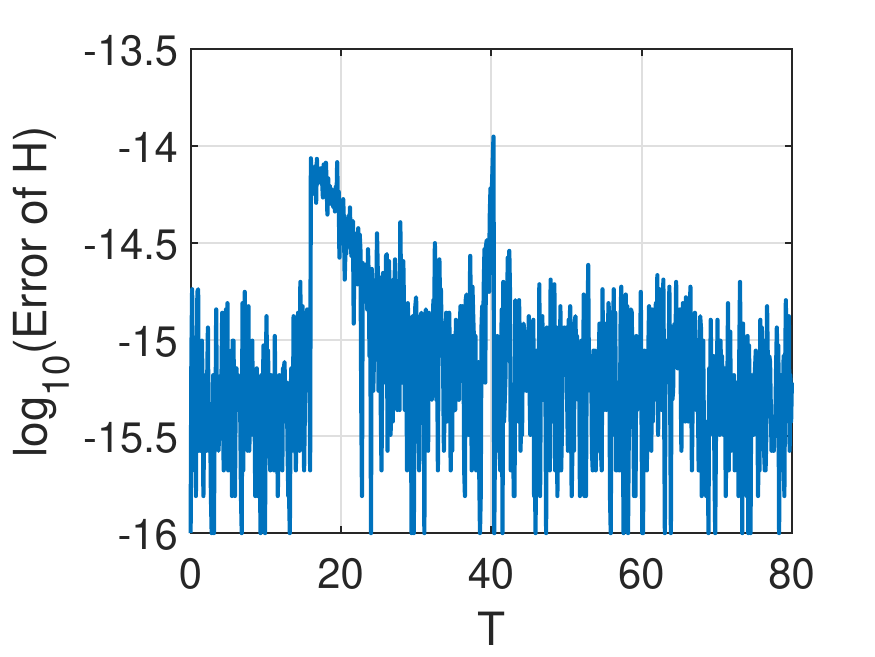}
		\end{minipage}
		\begin{minipage}[t]{0.32\textwidth}
			\centering
			\includegraphics[width=0.96\linewidth]{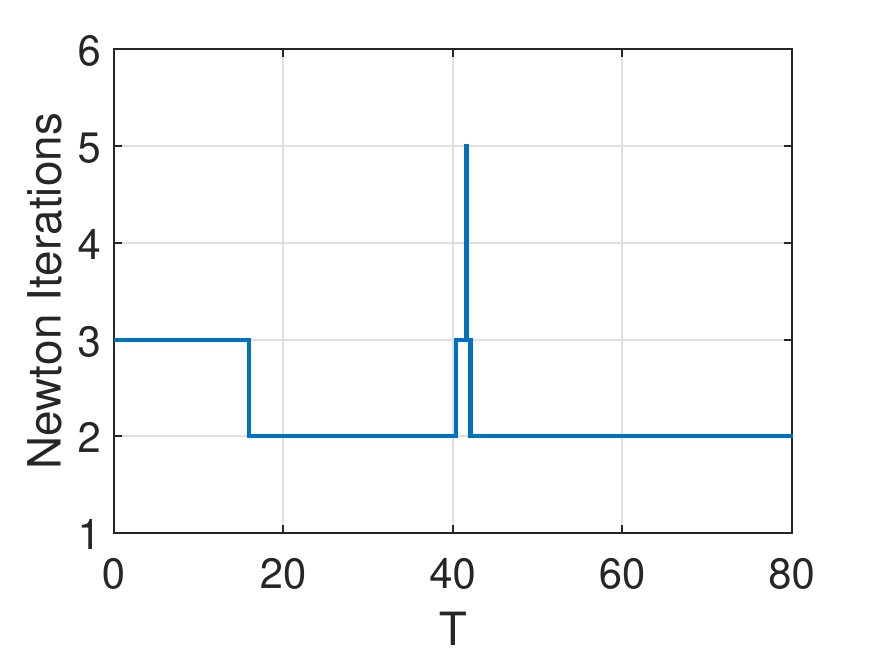}
		\end{minipage}
		\caption{Discrete energy (left), diffuse-volume conservation error (middle), and Newton iteration counts (right) for the migration and coalescence test.}
		\label{fig:migration-structure}
	\end{figure}

	Two droplets migrate toward a central hydrophilic region ($\theta_Y=20^\circ$) from outer regions with $\theta_Y=90^\circ$ and coalesce. Figures~\ref{selfclean}--\ref{fig:migration-structure} report shapes, contact histories, structural diagnostics, and Newton iterations.
	
	The computational domain is $\Omega=[-1,1]\times[0,0.5]$, with a uniform grid of $N_x=256$ and $N_y=64$, giving $\Delta x=\Delta y=1/128$. We take $\eps=0.02$, $\tau=0.01$, and the final time $T=80$. Periodic boundary conditions are imposed in the $x$-direction, the dynamic boundary condition is imposed on the lower wall, and the homogeneous Neumann condition is imposed at the upper boundary. The central hydrophilic region is defined by $|x|\le a$, with $a=0.6$, while the outer regions have $|x|>a$.
	
	Initially, two mirror-symmetric circular droplets of radius $R_0=0.13$ are centered at $x=\pm a$, with initial contact angle $\theta_0=110^\circ$. Thus, the initial droplets are located at the interfaces between the two substrate materials. The evolution is recorded at $t=0$, $16$, $40$, and $80$.
	
	Figure~\ref{selfclean} shows inward migration and coalescence near $t=41$, followed by relaxation of the merged droplet. On the left, Figure~\ref{fig:migration-contact} shows the outer contact approaching $x=-a$ and the inner contact moving toward $x=0$ before disappearing at coalescence; the right side is its mirror image. Both angles decrease from $110^\circ$, and the surviving apparent outer angle approaches $20^\circ$. Energy decreases, particularly around coalescence, while the diffuse-volume error stays near $10^{-14}$--$10^{-16}$ (Figure~\ref{fig:migration-structure}), demonstrating structure preservation through this topological change.
	
	
	\subsection{Wetting on rough substrates}
	
\begin{figure}[t!]
		\centering
		\begin{minipage}[t]{0.99\textwidth}
			\centering
			\includegraphics[width=0.96\textwidth,trim=0bp 25bp 0bp 25bp,clip]{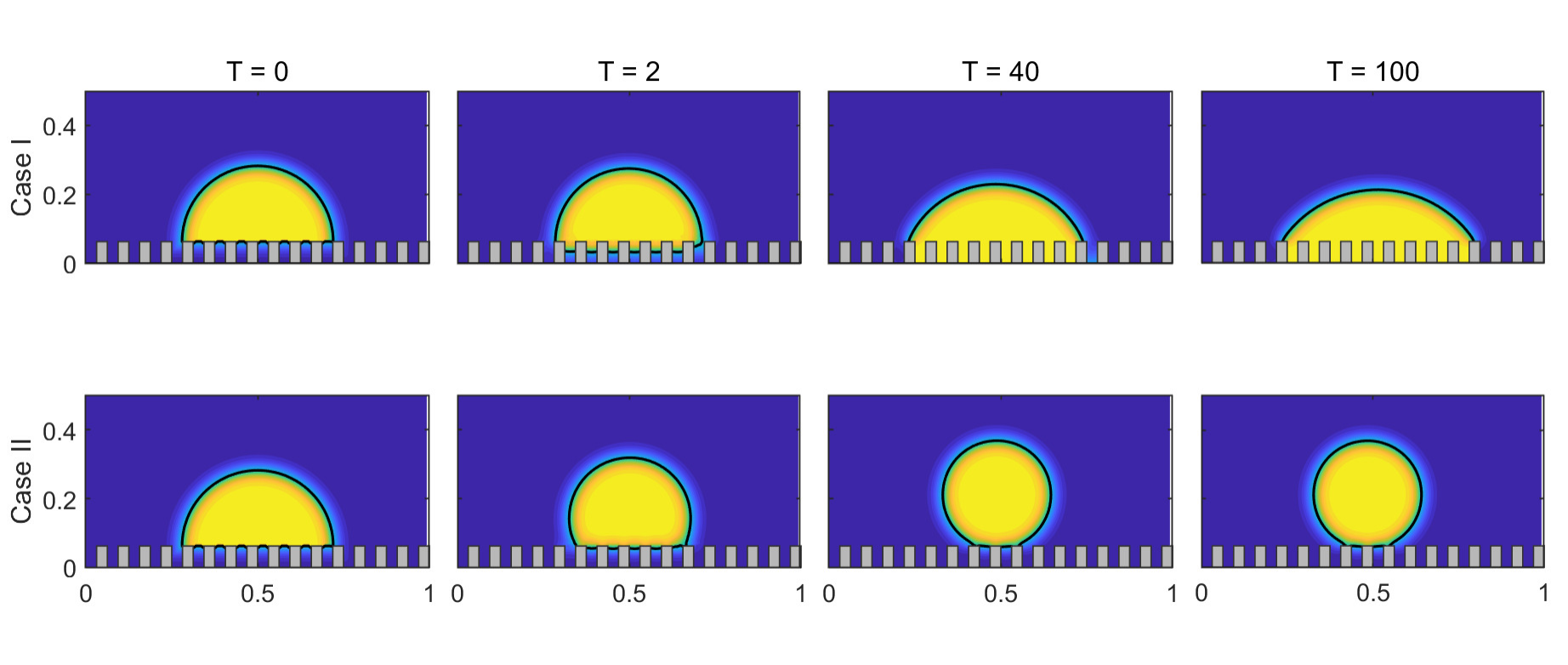}
		\end{minipage}
		\caption{Droplet evolution toward Wenzel-type (Case I) and Cassie--Baxter-type (Case II) configurations on rough substrates.}
		\label{fig:rough-evolution}
	\end{figure}

\begin{figure}[htbp]
		\centering
		\begin{minipage}[t]{0.32\textwidth}
			\centering
			\includegraphics[width=0.96\linewidth]{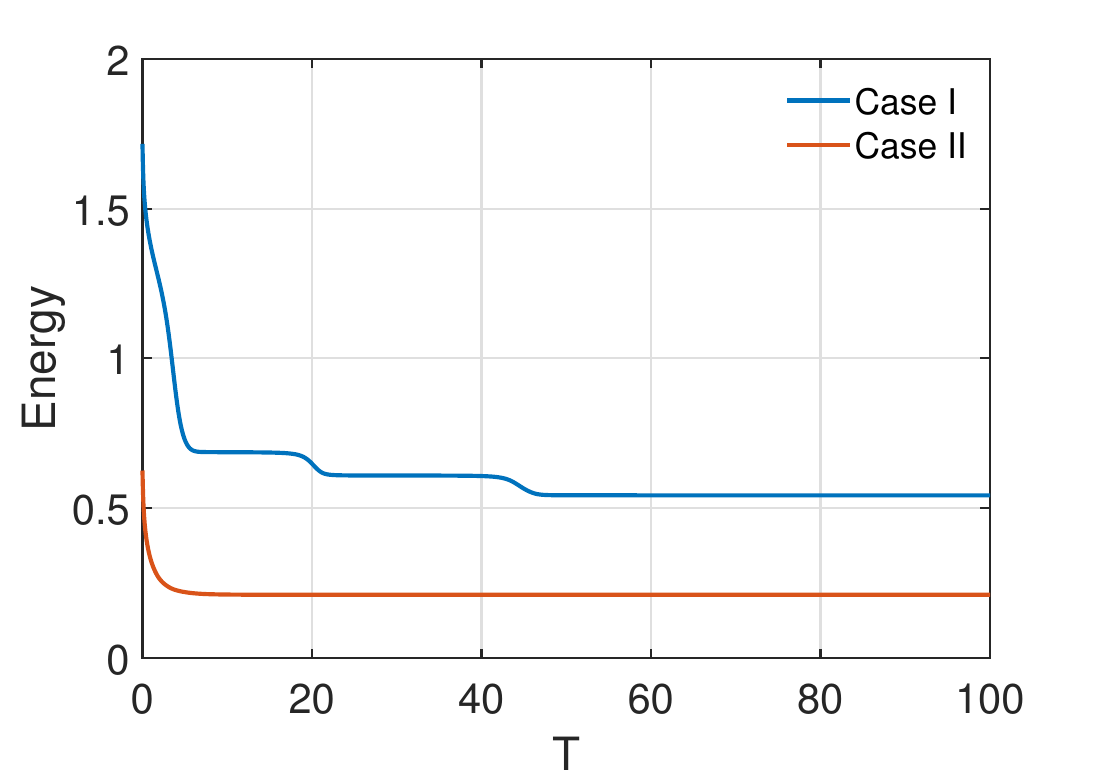}
		\end{minipage}
		\begin{minipage}[t]{0.32\textwidth}
			\centering
			\includegraphics[width=0.96\linewidth]{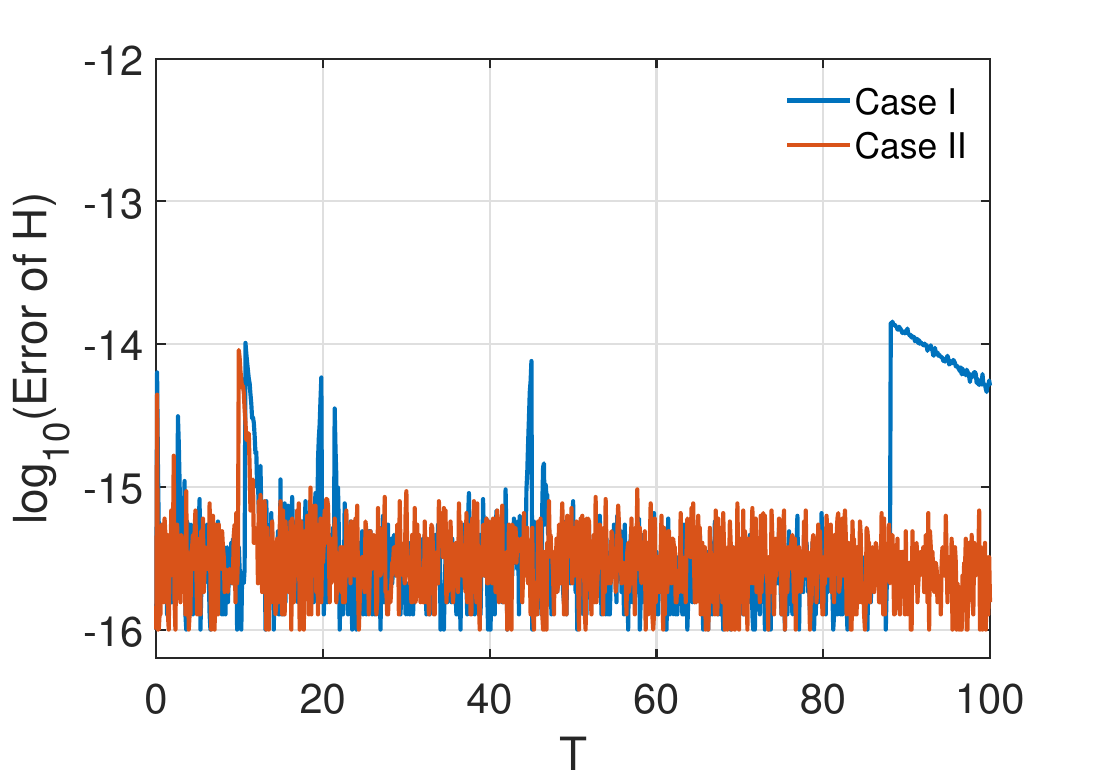}
		\end{minipage}
		\begin{minipage}[t]{0.32\textwidth}
			\centering
			\includegraphics[width=0.96\linewidth]{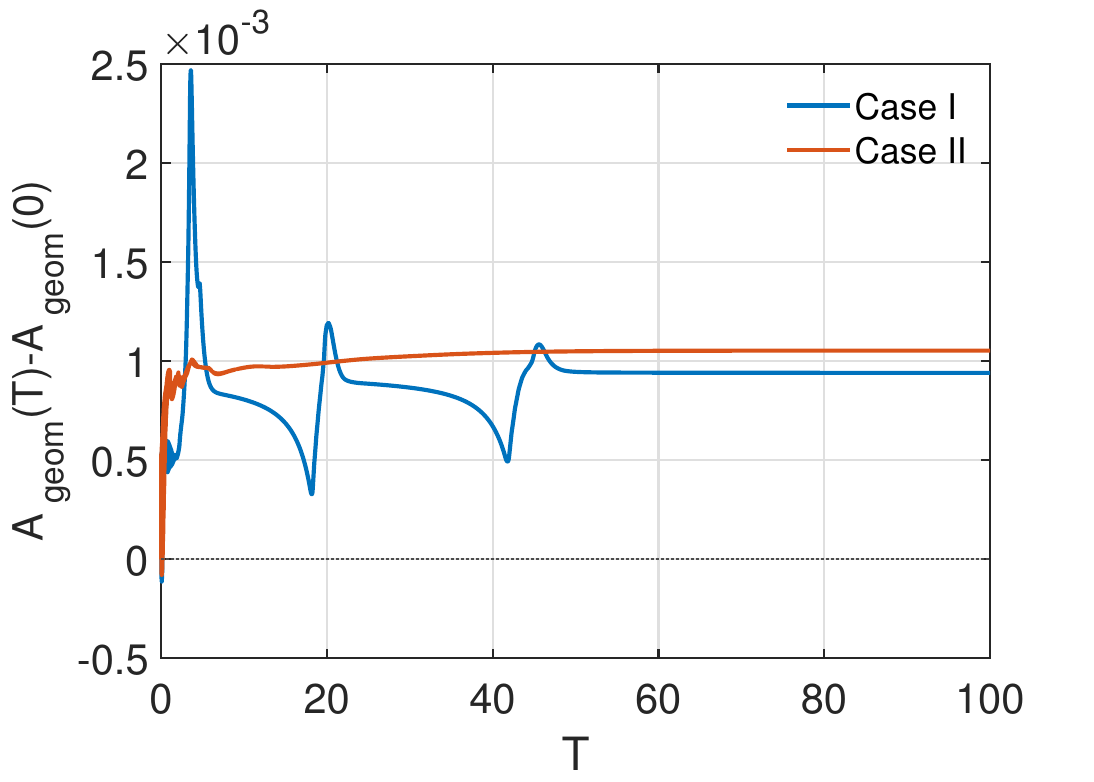}
		\end{minipage}
		\caption{Discrete energy (left), diffuse-volume conservation error (middle) and relative geometric-area error (right) for wetting on rough substrates.}
		\label{fig:rough-structure}
	\end{figure}

	We finally test periodic grid-aligned pillars, comparing penetration into grooves (Wenzel-type wetting) with retained vapor pockets (Cassie--Baxter-type wetting); see Figures~\ref{fig:rough-evolution} and~\ref{fig:rough-structure}.
	
	The computational domain is $\Omega=[0,1]\times[0,0.5]$, with $N_x=128$ and $N_y=64$, giving $\Delta x=\Delta y=1/128$. We take $\eps=0.02$, $\tau=0.01$, and the final time $T=100$. The pillars are aligned with the grid, with a width of $4\Delta x$, a gap of $4\Delta x$, and a height of $8\Delta y$. Periodic boundary conditions are imposed in the $x$-direction, the homogeneous Neumann condition is imposed at the upper boundary, and the dynamic wetting condition is imposed on the exposed substrate and pillar surfaces.
	
	The initial droplet is centered at $x=0.5$, with its nominal circular radius $R_0=0.22$ and its center at the pillar-top height. The grooves are initially filled with vapor, giving a Cassie-type initial configuration. We consider Young contact angles $\theta_Y=60^\circ$ and $120^\circ$, corresponding to Case I and Case II, respectively.
	
	Solid cells are excluded. Bulk weights come from fluid cells, boundary weights from exposed fluid--solid edges, and the gradient-energy matrix from fluid-domain grid edges. Thus, we construct the energy and its gradient consistently on the fluid region before applying the same discrete-variational-derivative scheme; the rectangular summation-by-parts identity is not used unchanged.
	
	We also measure $A_h(t)=|\{U>0\}|$ in the fluid domain. Each fluid cell is split into two triangles; piecewise-linear nodal interpolation gives their liquid fractions, whose sum yields $A_h$. The diagnostic $A_h(t)-A_h(0)$ is independent of the diffuse-volume constraint. Figure~\ref{fig:rough-evolution} shows groove penetration for $60^\circ$ (Case I) and retained vapor pockets for $120^\circ$ (Case II). Figure~\ref{fig:rough-structure} shows energy decay, small diffuse-volume error, and $O(10^{-3})$ changes in geometric area, which is not a discrete invariant.
	
	\FloatBarrier
	\section{Conclusions}\label{sec:conclusion}
   We developed a structure-preserving Allen–Cahn framework for three-phase contact with dynamic wall relaxation. The nonlinear diffuse-volume density varies quadratically near the pure-phase values, reducing the effect of bulk phase shifts on the measured volume. In the Onsager formulation, the volume gradient weights the correction toward the interface, and volume conservation makes its work vanish. The CVC scheme reproduces this cancellation through the discrete chain rule. Summation by parts couples the bulk and wall equations. As a result, the scheme exactly conserves the nonlinear diffuse volume, and the decrease in the original discrete energy equals the sum of bulk and wall dissipation at each step. Under the stated assumptions, the phase field in the continuous model remains within the physical interval, and coercivity of the discrete energy yields a mesh-dependent maximum-norm bound.
    
    The numerical experiments show second-order convergence in space and time and verify both discrete identities to solver accuracy. In the flat-wall comparison, the smallest droplet disappears under the linear mass constraint even though the absolute mass error remains near ($10^{-13}$) or below. The proposed formulation retains the droplet and yields an equilibrium shape close to the reference circular cap. Further tests capture contact-line pinning and wettability-driven migration and coalescence on chemically heterogeneous walls, as well as Wenzel- and Cassie–Baxter-type configurations on rough substrates.

	\bibliographystyle{siamunsrt}
	\begingroup
	\interlinepenalty=10000
	\bibliography{Bib_AC}
	\endgroup
	
	\appendix
	\section{Geometric and energetic interpretation of the nonlinear volume constraint}
	\label{app:diffuse-volume}
	
	Define the normalized density
	\begin{align*}
		h(s)=\frac{H(s)-H(-1)}{H(1)-H(-1)}.
	\end{align*}
	Then $h(-1)=0$, $h(1)=1$, and $V_H(u)=\int_\Omega h(u)\,\mathrm{d}x$. In particular, for $u=2\mathbf1_D-1$,
	\begin{align*}
		\mathcal V(u)=H(-1)|\Omega|+\bigl(H(1)-H(-1)\bigr)|D|,
		\qquad
		V_H(u)=|D|.
	\end{align*}
	Thus, $h(u)$ provides a diffuse approximation of the phase indicator, and its integral recovers the corresponding phase volume in the sharp-interface description.
	
	The choice $H'=\sqrt{2W}$ also makes this density compatible with the prescribed interfacial energy. While the integral of a phase indicator gives its volume, its total variation gives its perimeter. The corresponding connection for the diffuse density follows from the chain rule and the identity
	\begin{align*}
		\frac{\varepsilon}{2}|\nabla u|^2+\frac{W(u)}{\varepsilon}
		&=
		|\nabla H(u)|
		+\frac{1}{2\varepsilon}
		\left(\varepsilon|\nabla u|-\sqrt{2W(u)}\right)^2\\
		&\geq
		\bigl(H(1)-H(-1)\bigr)|\nabla h(u)|.
	\end{align*}
	Hence the bulk energy controls the total variation of the same normalized density used to measure phase volume. More precisely, if $u_\varepsilon\in H^1(\Omega)$ takes values in $[-1,1]$ and converges to $2\mathbf1_D-1$ in $L^1(\Omega)$, then we have $h(u_\varepsilon)\to\mathbf1_D$ in $L^1(\Omega)$. By the lower semicontinuity of total variation, one infers
	\begin{align*}
		\liminf_{\varepsilon\to0}
		\int_\Omega
		\left(
		\frac{\varepsilon}{2}|\nabla u_\varepsilon|^2
		+\frac{W(u_\varepsilon)}{\varepsilon}
		\right)\,\mathrm{d}x
		\geq
		\bigl(H(1)-H(-1)\bigr)\operatorname{Per}(D;\Omega).
	\end{align*}
	This relates the diffuse energy to the interfacial area in the sharp-interface limit.
	
	The square term above measures the departure from equipartition between the gradient and potential energies. It vanishes for an optimal one-dimensional transition. Moreover preciselt, the increasing profile $q$ satisfying $q''=W'(q)$ and $q(\pm\infty)=\pm1$ obeys
	\begin{align*}
		q'=\sqrt{2W(q)}=H'(q).
	\end{align*}
	For $u(x)=q(x_d/\varepsilon)$, the energy per unit interfacial area therefore equals
	\begin{align*}
		\int_{\mathbb R}
		\left(\frac12|q'|^2+W(q)\right)\,\mathrm{d}z
		=
		\int_{-1}^{1}\sqrt{2W(s)}\,\mathrm{d}s
		=
		H(1)-H(-1).
	\end{align*}
	Thus, the choice of $H$ is consistent with both the optimal transition profile and the surface tension associated with the energy; cf. \cite{ModicaMortola}. Along this profile, $H'(q)=q'$ also identifies the constraint gradient with the translational mode of the interface.
	
	Finally, the nonlinear density reduces the contribution of small bulk deviations to the measured phase volume. Assume that $W$ is $C^2$ near $[-1,1]$, with $W'(\pm1)=0$ and $W''(\pm1)>0$. It holds from the Taylor expansion that
	\begin{alignat*}{3}
		h(s)&=O((s+1)^2) &&\qquad\text{as }s\to-1^+,\\
		1-h(s)&=O((1-s)^2) &&\qquad\text{as }s\to1^-.
	\end{alignat*}
	In contrast, the linear density $(1+s)/2$ departs from its pure-phase values at first order. Since $W$ vanishes quadratically at the wells and is positive on $(-1,1)$, there exists a constant $C>0$, depending only on $W$, such that
	\begin{align*}
		\left|h(s)-\mathbf1_{\{s>0\}}\right|
		\leq C W(s),
		\qquad s\in[-1,1].
	\end{align*}
	Consequently, whenever $-1\leq u\leq1$ almost everywhere, we derive
	\begin{align*}
		\left|V_H(u)-|\{u>0\}|\right|
		&\leq
		\int_\Omega
		\left|h(u)-\mathbf1_{\{u>0\}}\right|\,\mathrm{d}x\leq C\int_\Omega W(u)\,\mathrm{d}x.
	\end{align*}
	With a uniform bound on the bulk diffuse-interface energy therefore, one obtains an $O(\varepsilon)$ bound on this discrepancy. Exact conservation of $V_H$ nevertheless remains distinct from exact preservation of the geometric volume $|\{u>0\}|$ at finite $\varepsilon$.
	
	\section{MBP continuous}
	\begin{proposition}\label{prop:continuous-interval}
		Let $u_\eps\in C^{2,1}(\overline\Omega\times[0,T])$ solve \eqref{eq:constrained-allen-cahn} with the prescribed periodic condition on the remaining boundary. Assume that $\lambda_\eps\in L^1(0,T)$. If $-1\le u_\eps(0)\le 1$ on $\overline\Omega$, then
		\begin{align*}
			-1\le u_\eps(x,t)\le 1
			\qquad \text{for all } (x,t)\in\overline\Omega\times[0,T].
		\end{align*}
	\end{proposition}
	
	\begin{proof}
		For $\sigma\in\{-1,1\}$, we define $w_\sigma=(\sigma u_\eps-1)_+$ and
		\begin{align*}
			E_\sigma(t)=\eps\|w_\sigma(t)\|_{L^2(\Omega)}^2+\|w_\sigma(t)\|_{L^2(\Gamma)}^2.
		\end{align*}
		Testing the bulk equation by $\eps\sigma w_\sigma$ and the dynamic boundary equation by $\sigma w_\sigma$, we obtain
		\begin{align*}
			\frac12 E_\sigma'(t)+\eps\|\nabla w_\sigma\|_{L^2(\Omega)}^2
			+\frac{\sigma}{\eps}\int_\Omega W'(u_\eps)w_\sigma\dx
			&=\eps\sigma\lambda_\eps(t)\int_\Omega H'(u_\eps)w_\sigma\dx\\
			&\quad-\sigma\int_\Gamma g'(u_\eps)w_\sigma\dS.
		\end{align*}
        Here, the normal-flux terms on $\Gamma$ cancel, while the remaining boundary terms vanish by the prescribed boundary conditions. 
		
		For the prescribed double-well potential, we known $\sigma W'(s)\ge 0$ whenever $\sigma s\ge 1$. Thus the potential term on the left-hand side is nonnegative. By choosing $M\ge 1$ such that $|u_\eps|\le M$ on $\overline\Omega\times[0,T]$, on $[-M,M]$, it holds $H'$ and $g'$ are Lipschitz in the phase variable, with constants $L_H$ and $L_g$, respectively, and $L_g$ can be chosen uniformly in $x\in\Gamma$. Since $H'(\sigma)=0$ and $g'(x,\sigma)=0$, on $\{w_\sigma>0\}$ we have
		\begin{align*}
			|H'(u_\eps)|&\le L_H|u_\eps-\sigma|=L_Hw_\sigma
			\qquad \text{on }\Omega,\\
			|g'(u_\eps)|&\le L_g|u_\eps-\sigma|=L_gw_\sigma
			\qquad \text{on }\Gamma.
		\end{align*}
		Consequently,
		\begin{align*}
			\frac12 E_\sigma'(t)+\eps\|\nabla w_\sigma\|_{L^2(\Omega)}^2
			&\le \eps L_H|\lambda_\eps(t)|\|w_\sigma\|_{L^2(\Omega)}^2
			+L_g\|w_\sigma\|_{L^2(\Gamma)}^2\\
			&\le \bigl(L_H|\lambda_\eps(t)|+L_g\bigr)E_\sigma(t).
		\end{align*}
		Because $E_\sigma(0)=0$ and the coefficient on the right-hand side is integrable, Gronwall's inequality yields $E_\sigma(t)=0$ for all $t\in[0,T]$. Taking both signs and using continuity of $u_\eps$ shows the claim.
	\end{proof}
	
	\begin{remark}[Integrability of the quotient-defined multiplier]
		Suppose that $\Omega$ is connected, that $\lambda_\eps$ is given by \eqref{eq:lambda} whenever its denominator is positive, and that $u_\eps(0)$ takes values in $[-1,1]$ and is neither identically $1$ nor identically $-1$. Set $D(t)=\int_\Omega 2W(u_\eps(t))\dx$. Then $D(0)>0$ and
		\begin{align*}
			|\Omega|H(-1)<\cV(u_\eps(0))<|\Omega|H(1).
		\end{align*}
		If $D$ has a first zero at $t_*\in(0,T]$, the bulk equation and \eqref{eq:lambda} give, for $0<t<t_*$,
		\begin{align*}
			\frac{\mathrm d}{\mathrm dt}\cV(u_\eps(t))
			&=\int_\Omega H'(u_\eps)
			\bigl(\Delta u_\eps-\eps^{-2}W'(u_\eps)\bigr)\dx
			+\lambda_\eps(t)D(t)=0.
		\end{align*}
		By continuity, $\cV(u_\eps(t_*))=\cV(u_\eps(0))$. On the other hand, $D(t_*)=0$ implies that $u_\eps(\cdot,t_*)$ takes values only in $\{-1,1\}$. By continuity and connectedness, it is identically $1$ or $-1$, contradicting the strict bounds on $\cV(u_\eps(0))$. Hence $\min_{[0,T]}D>0$. The numerator in \eqref{eq:lambda} is continuous on $[0,T]$, so $\lambda_\eps\in L^\infty(0,T)\subset L^1(0,T)$. For the pure-phase initial states, $u_\eps\equiv\pm1$ are stationary solutions under the convention $\lambda_\eps=0$.
	\end{remark}
\end{document}